\documentclass[preprint,12pt]{elsarticle}
\usepackage{times}
\usepackage{amsfonts,amsmath,amssymb,amsthm}
\usepackage{color,soul}
\biboptions{numbers,sort&compress} 
\let\ts=\textstyle
\def\fracskip{\mskip 1mu \relax}
\let\oldfrac=\frac
\def\nfrac#1#2{\oldfrac{\fracskip#1\fracskip}{\fracskip#2\fracskip}}
\def\tfrac#1#2{{\ts\nfrac{#1}{#2}}}
\let\frac=\nfrac
\def\pd#1#2{\frac{\partial#1}{\partial#2}}
\def\pdd#1#2#3{\ifx#2#3\pd{^2#1}{#2^2}\else\pd{^2#1}{#2\partial#3}\fi }
\def\dash{\unskip\nobreak\hskip.05555em---\hskip.05555em\relax}

\let\Bl=\Bigl \let\Br=\Bigr
\let\BL=\biggl \let\BR=\biggr
\let\bl=\bigl \let\br=\bigr

\def\arbs{are arbitrary constants}

\def\Equation#1. {\medbreak{\bfseries\itshape{Equation\kern.3333em\relax#1.}}\enspace\ignorespaces }

\newtheoremstyle{remark}{\medskipamount}{\medskipamount}
  {\small\rmfamily}{\parindent}{\footnotesize\sffamily}{.}{.5em}{}
\theoremstyle{remark}
\newtheorem{remark}{Remark}

\newtheoremstyle{example}{\medskipamount}{\medskipamount}
  {\small\rmfamily}{\parindent}{\footnotesize\sffamily}{.}{.5em}{}
\theoremstyle{example}
\newtheorem{example}{Example}

\let\ds=\displaystyle
\let\ts=\textstyle
\let\bl=\bigl \let\br=\bigr
\let\Bl=\Bigl \let\Br=\Bigr
\let\BL=\biggl \let\BR=\biggr

\begin{document}
\begin{frontmatter}
\title{Numerical integration of blow-up problems\\
on the basis of non-local transformations\\
and differential constraints$^*$}
\author[ipm]{Andrei D. Polyanin\corref{cor1}}
\author[us]{Inna K. Shingareva\corref{cor2}}
 \cortext[]{This article is a significantly extended and more general version of the article A.D.~Polyanin and I.K. Shingareva,
The use of differential and non-local transformations for numerical integration of non-linear blow-up problems,
\textit{Int. J. Non-Linear Mechanics}, 2017, Vol.~95, pp.~178--184.}
\address[ipm]{Institute for Problems in Mechanics, Russian Academy of Sciences,\\
  101 Vernadsky Avenue, bldg~1, 119526 Moscow, Russia}
\address[bmstu]{Bauman Moscow State Technical University,\\
5 Second Baumanskaya Street, 105005 Moscow, Russia}
\address[mephi]{National Research Nuclear University MEPhI,
  31 Kashirskoe Shosse, 115409 Moscow, Russia}
\address[us]{University of Sonora, Blvd. Luis Encinas y Rosales S/N, Hermosillo C.P. 83000, Sonora, M\'exico}

\begin{abstract}
Several new methods of numerical integration of Cauchy problems with blow-up solutions
for nonlinear ordinary differential equations of the first- and second-order are described.
Solutions of such problems have singularities whose positions are unknown a priori
(for this reason, the standard numerical methods for solving problems with blow-up solutions
can lead to significant errors). The first proposed method is based on the
transition to an equivalent system of equations by introducing
a new independent variable chosen as the first derivative,
$t=y^{\prime}_x$, where $x$ and $y$ are independent and dependent variables in the original equation.
The second method is based on introducing a new auxiliary non-local variable
of the form  $\xi=\int^x_{x_0} g(x,y,y^{\prime}_x)\,dx$  with the subsequent transformation to
the Cauchy problem for the corresponding system of ODEs.
The third method is based on adding to the original equation of a differential constraint,
which is an auxiliary ODE connecting the given variables and a new variable.
The proposed methods lead to problems whose solutions are represented in parametric form and
do not have blowing-up singular points;
therefore the transformed problems admit the application of standard fixed-step numerical methods.
The efficiency of these methods is illustrated by solving a number of test problems that admit
an exact analytical solution.
It is shown that: (i) the methods based on non-local transformations of a special kind
are more efficient than several other methods, namely, the method based on the hodograph transformation,
the method of the arc-length transformation, and the method based on the differential transformation, and
(ii) among the proposed methods, the most general method is the method based on the differential constraints.
Some examples of nonclassical blow-up problems are considered, in which the right-hand side
of equations has fixed singular points or zeros.
Simple theoretical estimates are derived for the critical value of an independent variable
bounding the domain of existence of the solution.
It is shown that the method based on a non-local transformation of the general form
as well as the method based on the differential constraints admit generalizations to the
$n$\,th-order ODEs and systems of coupled ODEs.
\end{abstract}

\end{frontmatter}

\section{Introduction}\label{ss:1}

\subsection{Preliminary remarks. Blow-up solutions}\label{ss:1.1}

We will consider Cauchy problems for ordinary differential equations (briefly, ODEs), whose solutions
tend to infinity at some finite value of the independent variable $x=x_*$, where $x_*$
does not appear explicitly in the differential equation under consideration and it is not known in advance.
Similar solutions exist on a bounded interval (hereinafter in this article we assume that $x_0\le x<x_*$)
and are called blow-up solutions.
This raises the important question for practice: how to determine the position of a singular point~$x_*$
and the solution in its neighborhood using numerical methods.

In the general case, the blow-up solutions that have a power singularity
can be represented in a neighborhood  of the singular point~$x_*$ in the form
\begin{align}
y\approx A(x_*-x)^{-\beta},\quad \ \beta>0,
\label{eq:01}
\end{align}
where $A$ and $\beta$ are some constants.
For these solutions we have $\ds\lim_{x\to x_*}y=\infty$ and $\ds\lim_{x\to x_*}y^{\prime}_x=\infty$.

For blow-up solutions with the power singularity~\eqref{eq:01} near the singular point~$x_*$ we have
\begin{align}
\frac{y^{\prime}_x}{y} \approx \frac{\beta}{x_*-x},
\label{eq:01-2}
\end{align}
i.e. the required function grows more slowly than its derivative. Therefore,
$\ds\lim_{x\to x_*} |y^{\prime}_x/y|=\infty$ (this is a common property of any blow-up solutions).

\begin{example}
Consider the test Cauchy problem for the first-order nonlinear ODE
with separable variables
\begin{align}
y^{\prime}_x=y^2\quad (x>0),\qquad y(0)=1.
\label{eq:02f}
\end{align}
The exact solution of this problem has the form
\begin{align}
y=\frac{1}{1-x}.
\label{eq:02g}
\end{align}
It has a power-type singularity (a first-order pole) at the point $x_*=1$
and does not exist for $x>x_*$.

The Cauchy problem \eqref{eq:02f} is a particular case of the three-parameters problem
\begin{align}
y^{\prime}_x=by^{\gamma}\quad (x>0),\quad \ y(0)=a,
\label{eq:xx02}
\end{align}
where $a$, $b$, and $\gamma$ \arbs.
If the inequalities
\begin{align}
a>0,\quad \ b>0, \quad \ \gamma>1
\label{eq:xx03}
\end{align}
are valid, then the exact solution of the problem~\eqref{eq:xx02} is given by the formula
\begin{align}
y=A(x_*-x)^{-\beta},
\label{eq:xx04}
\end{align}
where
$$
A=[b(\gamma-1)]^{\ts\frac 1{1-\gamma}},\quad \ x_*=\frac 1{a^{\gamma-1}b(\gamma-1)},\quad \ \beta=\frac 1{\gamma-1}>0.
$$
This solution exists on a bounded interval $0\le x<x_*$, where $x_*$ is a singular point
of the pole-type solution, and does not exist for $x\ge x_*$.
In this case, the solution~\eqref{eq:xx04} coincides with its asymptotic behavior in a neighborhood of the singular point
(compare \eqref{eq:01} with~\eqref{eq:xx04}).
\end{example}

There exist problems that have blow-up solutions with a different type of singularity (that differs from \eqref{eq:01}).
In particular, solutions with a logarithmic singularity at the point~$x_*$ have the form
\begin{align}
y\approx A\ln\bl[B(x_*-x)\br],
\notag
\end{align}
where $A$ and $B>0$ are some constants.

\begin{example}
The test Cauchy problem with exponential nonlinearity
\begin{align}
y^{\prime}_x=be^y\quad (x>0),\qquad y(0)=a
\label{*12.1.10.33}
\end{align}
admits the exact solution with a logarithmic singularity
\begin{align}
y=-\ln(e^{-a}-bx)
\label{*12.1.10.34}
\end{align}
for $a\ge 0$ and $b>0$.
This solution exists on the interval $0\le x<x_*=e^{-a}/b$ and does not exist for
$x\ge x_*$.
\end{example}

\subsection{Problems arising in numerical solutions of blow-up problems}\label{ss:1.2}

The direct application of the standard fixed-step numerical methods to blow-up problems
leads to certain difficulties because their solutions have a singularity
and the range of variation of the independent variable is unknown in advance~\cite{stuart1990}.
The difficulties arising in the application of the classical Runge--Kutta methods
for solving the test Cauchy problem~\eqref{eq:xx02} are described below
(the results of~\cite{als2005} are used).

The qualitative behavior of the numerical blow-up solution for equations of the form~\eqref{eq:xx02}
for $a>0$, $b>0$, and $\gamma>1$ is significantly different for the explicit and implicit Runge--Kutta methods
(the explicit methods up to the fourth-order of approximation and the Euler implicit
method have been tested in~\cite{als2005}).

All the explicit methods provide monotonically increasing solutions; and
the higher order of the approximation method, the
faster growth of the numerical solution. Soon after passing through the singular
point~$x_*$, in which the exact solution has a pole, an overflow occurs in the calculation
and further computing is impossible.
Such qualitative behavior is unpleasant, since it is difficult for the researcher to determine the cause of the overflow.

For the implicit methods, the picture is different. First, the solution increases, but
even before the pole it breaks down to the region of negative values.
The calculation of the right-hand side of~\eqref{eq:xx02} for fractional values of~$\gamma$
becomes impossible (because a fractional power of a negative number occurs).

Various special methods have been proposed in the literature
for numerical integration of problems that have blow-up solutions.

One of the basic ideas of numerical integration of blow-up problems
consists in the application of an appropriate transformation at the initial stage,
which leads to the equivalent problem for one differential equation or a system
of coupled equations whose solutions have no singularities at a priori unknown point
(after such transformations, the unknown singular point $x=x_*$ usually goes to
the infinitely remote point for the new independent variable).

Currently, two methods based on this idea are most commonly used.
The first method, based on the hodograph transformation, $x=\bar y$, $y=\bar x$
(where the independent and dependent variables are interchanged) was proposed in~\cite{aco2002}.
The second method of this kind, called the method of the arc-length transformation,
is described in~\cite{mor1979} (for details, see below Item~$2^\circ$ in Sections~3.1 and~7.1,
as well as reference~\cite{hir2006}).
This method is rather general and it can be applied for numerical integration of systems
of ordinary differential equations.

The methods based on the hodograph and arc-length transformations for blow-up
solutions with a power singularity of the form~\eqref{eq:01} lead to the
Cauchy problems whose solutions tends to the asymptote with respect to the
power law for large values of the new independent variable. This creates
certain difficulties in some problems, since one has to consider large
intervals of variation of the independent variable in numerical integration.

Based on other ideas, some special methods of numerical integration of
blow-up problems are described, for example, in
\cite{stuart1990,als2005,hir2006,meyer98,gor1999,eli2006,bar2008,nas2009,dlam2012,zhou2016,tak2017}.
In particular, it was suggested in \cite{eli2006,tak2017} to investigate such problems via compactifications,
which are point transformations of the special form (whose inverse transformations have singularities).

In this paper, we propose several new methods of numerical integration of Cauchy problems
for the first- and second-order nonlinear equations, which have blow-up solutions.
These methods are based on differential and non-local transformations, and also on differential constraints,
allowing us to obtain the equivalent problems for systems of equations whose solutions
do not have singularities at a priori unknown point.
Some special methods based on non-local transformations and differential constraints lead
to the Cauchy problems whose solutions, which are found in parametric form by numerical integration,
tend exponentially to the asymptote for large values of the new independent variable.
Therefore, these methods are more effective than
the methods based on the hodograph and arc-length transformations,
which lead to solutions that are quite slowly (by the power law) tend to the asymptote
for large values of the independent variable.
The presentation of the material is widely illustrated with test problems that admit an exact solution.
Two-sided theoretical estimates are established for the critical value of the independent variable $x=x_*$,
when an unlimited growth of the solution occurs as approaching it.
It is shown that the method based on a non-local transformation of the general form
as well as the method based on the differential constraints admit generalizations to
the $n$ th-order ordinary differential equations and systems of differential equations.

\section{Problems for first-order equations. Differential transformations}\label{ss:2}

\subsection{Solution method based on introducing a differential variable}\label{ss:2.1}

The Cauchy problem for the first-order differential equation has the form
\begin{align}
&y^{\prime}_x=f(x,y)\quad (x>x_0),\label{eq:02}\\
&
y(x_0)=y_0.\label{eq:02a}
\end{align}
In what follows we assume that $f=f(x,y)>0$, $x_0\ge 0$, $y_0>0$, and also $f/y\to\infty$ as $y\to\infty$
(in such problems, blow-up solutions arise when the right-hand side of a nonlinear equation
is quite rapidly growing as $y\to\infty$).

First, we represent the nonlinear ODE~\eqref{eq:02} in the form of an equivalent system
of differential-algebraic equations
\begin{align}
t=f(x,y),\quad \ y^{\prime}_x=t,
\label{eq:02b}
\end{align}
where $y=y(x)$ and $t=t(x)$   are unknown functions to be determined.

By applying~\eqref{eq:02b} and assuming that $y=y(t)$ and $x=x(t)$,
we derive a system of ODEs of the standard form.
By taking the full differential of the first equation
of~\eqref{eq:02b} and multiplying the second equation by~$dx$, we get
\begin{equation}
dt=f_x\,dx+f_y\,dy,\quad dy=t\,dx,
\label{eq:02c}
\end{equation}
where $f_x$ and $f_y$ denote the corresponding partial derivatives of the function $f=f(x,y)$.
Eliminating first~$dy$ and then~$dx$ from~\eqref{eq:02c}, we arrive at
the ODE system of the first order
\begin{equation}
x^{\prime}_t=\frac 1{f_x+tf_y},\quad \ y^{\prime}_t=\frac t{f_x+tf_y}\quad \ (t>t_0),
\label{eq:02d}
\end{equation}
which must be supplemented by the initial conditions
\begin{align}
x(t_0)=x_0,\quad y(t_0)=y_0,\quad t_0=f(x_0,y_0).
\label{eq:02e}
\end{align}
Conditions~\eqref{eq:02e} are derived from~\eqref{eq:02a} and the first equation of~\eqref{eq:02b}.

Assuming that the conditions $f_x+tf_y>0$ are valid at $t_0<t<\infty$, the Cauchy problem
\eqref{eq:02d}--\eqref{eq:02e} can be integrated numerically, for example, by applying
the Runge--Kutta method or other standard fixed-step numerical methods (see, for example, \cite{but,fox87,dormand96,sch,shampine94,asch,korn,shing2009,grif2010}).
In this case, the difficulties (described in Section 1.2) will not occur because of
the presence of a singularity in the solutions (since $x^{\prime}_t\to 0$ as $t\to\infty$).
In view of~\eqref{eq:02b}, the singular point~$x_*$ of the solution  corresponds to $t=\infty$,
therefore the required value~$x_*$ is determined by the asymptotic behavior of the function $x=x(t)$ for large~$t$.

\begin{remark}
Taking into account the first equation of~\eqref{eq:02b}, the system~\eqref{eq:02d}
can be represented in the form
\begin{equation}
x^{\prime}_t=\frac 1{f_x+tf_y},\quad \ y^{\prime}_t=\frac f{f_x+tf_y}\quad \ (t>t_0).
\label{eq:02dd}
\end{equation}
Another equivalent system of ODEs can be obtained by replacing $t$ by $f$ in~\eqref{eq:02dd}.
\end{remark}

\subsection{Test problem and numerical solutions}\label{ss:2.2}

Let us illustrate the method described in Section~2.1 with a simple example.

\begin{figure}
\centering
{\includegraphics[scale=0.34]{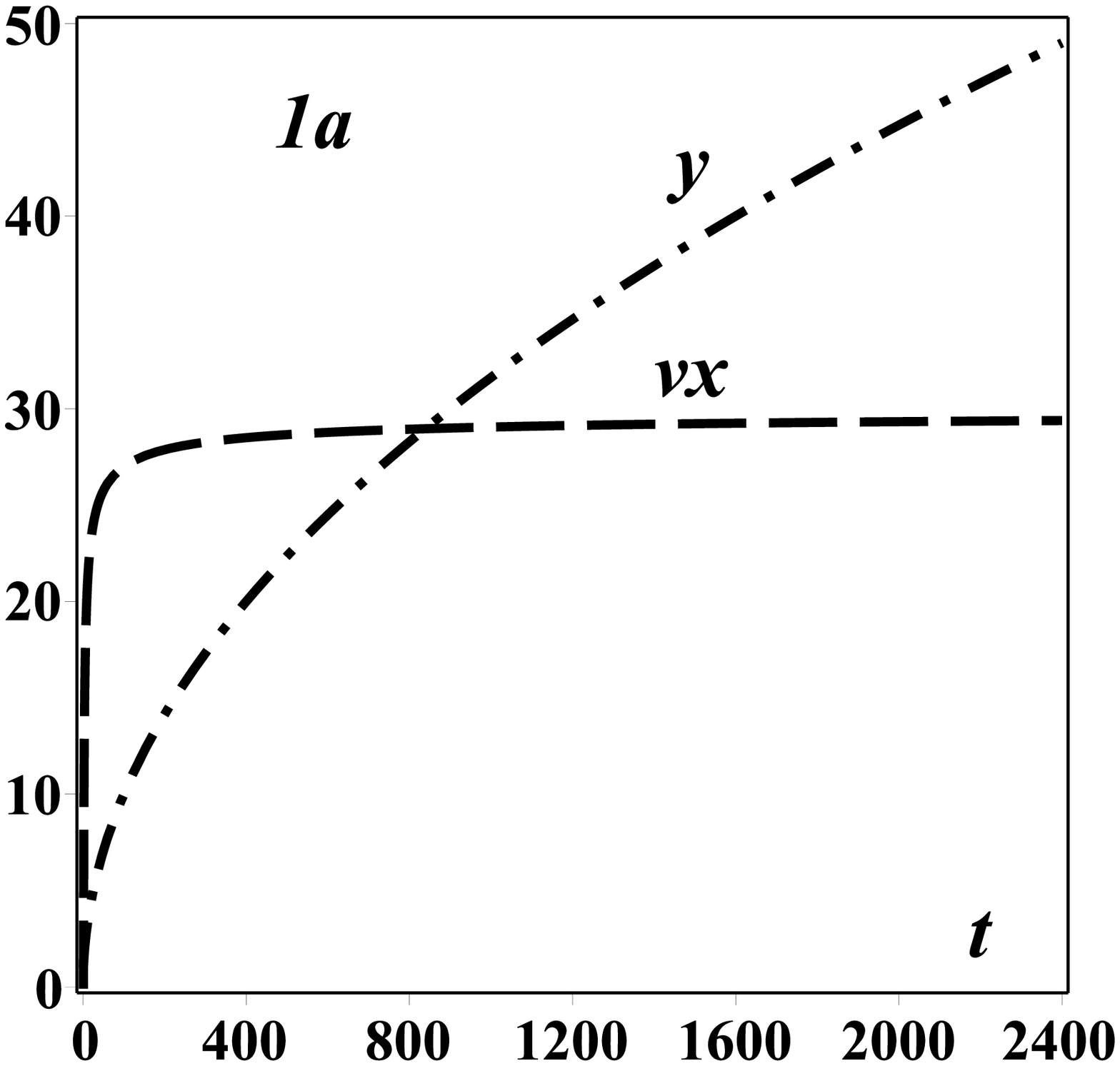}\ \includegraphics[scale=0.34]{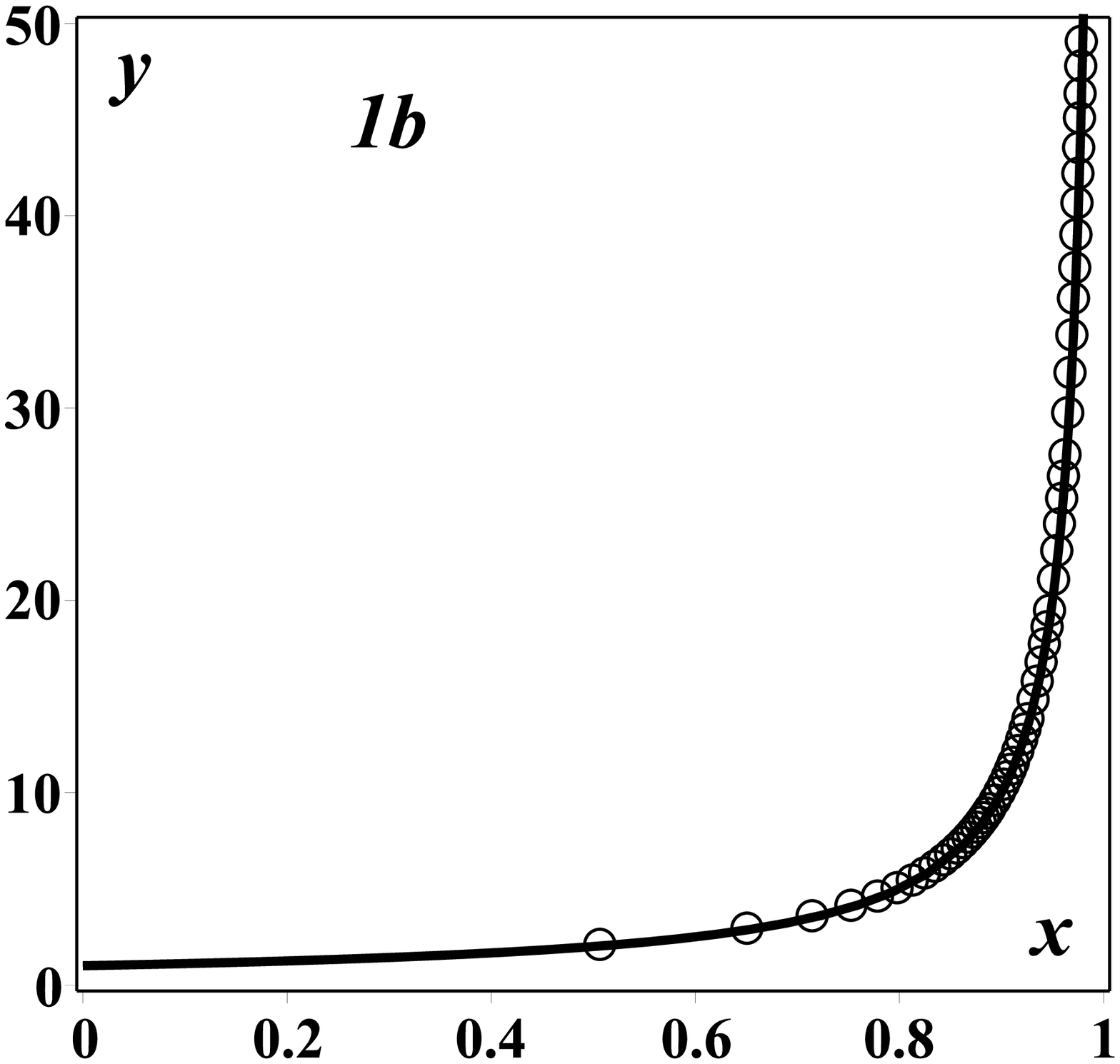}}
\caption{{\it 1a}\dash the dependences $x=x(t)$ and $y=y(t)$ obtained by numerical solution of the problem~\eqref{eq:02h}
for $a=b=1$ and $\gamma=2$ ($\nu=30$); {\it 1b}\dash exact solution~\eqref{eq:02g} (solid line) and numerical solution
of the problem~\eqref{eq:02h} for $a=b=1$ and $\gamma=2$ (circles).}
\label{fig:Fig1}
\end{figure}

\begin{example}
Consider the test Cauchy problem \eqref{eq:xx02}--\eqref{eq:xx03}.
By introducing a new variable $t=y^{\prime}_x$ in~\eqref{eq:xx02}, we obtain the following Cauchy
problem for the system of equations:
\begin{equation}
\begin{gathered}
x^{\prime}_t=\frac 1{b\gamma ty^{\gamma-1}},\quad \ y^{\prime}_t=\frac 1{b\gamma y^{\gamma-1}}\quad \ (t>t_0);\\
x(t_0)=0,\quad \ y(t_0)=a,\quad \ t_0=a^\gamma b,
\end{gathered}
\label{eq:02h}
\end{equation}
which is a particular case of the problem \eqref{eq:02d}--\eqref{eq:02e} for $f=by^\gamma$, $x_0=0$, and $y_0=a$.
The exact solution of the problem~\eqref{eq:02h} has the form
\begin{align}
x=\frac 1{b(\gamma-1)}\BL[a^{1-\gamma}-\Bl(\frac bt\Br)^{\!\ts\frac{\gamma-1}\gamma}\BR], \quad \ y=\Bl(\frac tb\Br)^{\ts\frac 1\gamma}\qquad (t\ge a^\gamma b).
\label{eq:02j}
\end{align}
It has no singularities; the function $x=x(t)$ increases monotonically with
$t>a^\gamma b$, tending to the desired limiting value $\ds x_*=\lim_{t\to\infty}x(t)=\frac 1{a^{\gamma-1}b(\gamma-1)}$,
and the function $y=y(t)$ monotonously increases with increasing~$t$.
The solution~\eqref{eq:02j} for the system~\eqref{eq:02h} is a solution of the original problem~\eqref{eq:xx02}--\eqref{eq:xx03} in
parametric form.

In Fig.~\ref{fig:Fig1}, we compare the exact solution~\eqref{eq:02g} of the Cauchy problem
for one equation~\eqref{eq:02f}  with the numerical solution of the transformed problem for
the system of equations~\eqref{eq:02h} for $a=b=1$ and $\gamma=2$, obtained by the classical numerical method, e.g. the
Runge--Kutta method of the fourth-order of approximation with a fixed step of integration,
equal to~$0.2$ (here and in what follows in the figures, for the sake of clarity,
a scale factor $\nu=30$ is introduced for the functions $x=x(t)$ or $x=x(\xi)$).
In this case, the maximum error of the numerical solution
does not exceed $0.017\%$ for $y\le 50$.
\end{example}

\begin{remark}
Here and in what follows, the numerical integration interval for the new
variable $t$ (or $\xi$) is usually determined, for demonstration
calculations, from the condition $\Lambda_\text{m}=50$, where

\begin{equation}
\Lambda_\text{m}=\min[y,\, y^{\prime}_x/y]\quad \ \text{(for $y_0\sim 1$ and $y_1=y'_x(x_0)\sim 1$)}.
\label{eq:Lambda}
\end{equation}
In a few cases, the
condition $\Lambda_\text{m}=100$ or $\Lambda_\text{m}=150$ is used, which is specially stipulated.
For first-order ODE problems of the form~\eqref{eq:02}--\eqref{eq:02a}, this definition of $\Lambda_\text{m}$ can be replaced
by the equivalent definition $\Lambda_\text{m}=\min[y,\, f/y]$.

Conditions $y_0\sim 1$ and $y_1\sim 1$ in \eqref{eq:Lambda} are not strongly essential, since the substitution $y=y_0-1+(y_1-1)(x-x_0)+\bar y$
leads to an equivalent problem with the initial conditions $\bar y(x_0)=\bar y^\prime_x(x_0)=1$.
\end{remark}

\subsection{Modified differential transformation}

The solution \eqref{eq:02j} tends rather slowly to the asymptotic values
$x\to x_*$ as $t\to \infty$ (in particular, for $s=2$ and large~$t$ we have $x_*-x\sim t^{-1/2}$).
To speed up the process of approaching the asymptotic behavior with respect to~$x$ in the system~\eqref{eq:02d}
is useful additionally to make the exponential-type substitution
\begin{align}
t=t_0\exp(\lambda\tau),\quad \ \tau\ge 0,
\label{eq:02k}
\end{align}
where
\begin{align}
\tau=\frac 1\lambda\ln\frac t{t_0}=\frac 1\lambda\ln\frac{y^{\prime}_x}{t_0}
\label{22222}
\end{align}
is a new independent variable and $\lambda>0$ is a numerical parameter that can be varied.
Transformations with a new independent variable of the form \eqref{22222} will be called the modified differential transformations.

\begin{example}
As a result of the substitution~\eqref{eq:02k}, the Cauchy problem~\eqref{eq:02h} is transformed to the form
\begin{equation}
\begin{gathered}
x^{\prime}_\tau=\frac\lambda{b\gamma y^{\gamma-1}},\quad \ y^{\prime}_\tau=\frac{a^\gamma\lambda e^{\lambda\tau}}{\gamma y^{\gamma-1}}\quad \ (\tau>0);\\
x(0)=0,\quad \ y(0)=a,
\end{gathered}
\label{eq:02l}
\end{equation}
and its exact solution is given by the formulas
\begin{align}
x=\frac {1}{a^{\gamma-1}b(\gamma-1)}\BL\{1-\exp\BL[-\frac{\lambda(\gamma-1)}\gamma \tau\BR]\BR\},\quad \ y=a\exp\BL(\frac{\lambda}\gamma\tau\BR),\qquad \tau\ge 0.
\label{eq:02m}
\end{align}

Let $a=b=1$, $\gamma=2$ and the stepsize is equal to~$0.4$.
For numerical integration of the test problem~\eqref{eq:02l} for $\lambda=1$ and $\lambda=2$
with the maximum error~$0.002\%$, it is required to take, respectively, the interval
$[0, 8]$ and $[0, 4]$ with respect to~$\tau$ to
approach the asymptote (however, for numerical integration of the related problem~\eqref{eq:02h}
with the maximum error~$0.016\%$, it is required to take an essentially larger interval $[0, 2980]$
with respect to~$t$).
\end{example}

\section{Problems for first-order equations. Nonlocal transformations and\\ differential constraints}\label{ss:3}

\subsection{Solution method based on introducing a non-local variable}\label{ss:3.1}

Introducing a new non-local variable~\cite{zai1993,zai1994,mur2010} according to the formula,
\begin{align}
\xi=\int^x_{x_0}g(x,y)\,dx,\quad \ y=y(x),
\label{*}
\end{align}
leads the Cauchy problem for one equation \eqref{eq:02}--\eqref{eq:02a} to the equivalent problem
for the autonomous system of equations
\begin{equation}
\begin{gathered}
x^{\prime}_\xi=\frac 1{g(x,y)},\quad \ y^{\prime}_\xi=\frac{f(x,y)}{g(x,y)}\quad \ (\xi>0);\\
x(0)=x_0,\quad \ y(0)=y_0.\qquad\qquad
\end{gathered}
\label{eq:02p}
\end{equation}
Here, the function $g=g(x,y)$ has to satisfy the following conditions:
\begin{equation}
g>0\ \text{for}\ x\ge x_0, \ y\ge y_0; \quad \ g\to \infty \ \text{as}\ y\to\infty;\quad \ f/g=k\ \text{as}\ y\to\infty,
\label{eq:02q}
\end{equation}
where $k=\text{const}>0$ (moreover, the limiting case $k=\infty$ is also allowed);
otherwise the function~$g$ can be chosen rather arbitrarily.

From~\eqref{*} and the second condition~\eqref{eq:02q} it follows that $x^{\prime}_\xi\to 0$ as $\xi\to\infty$.
The Cauchy problem~\eqref{eq:02p} can be integrated numerically applying the Runge--Kutta method
or other standard numerical methods.

\medskip

Let us consider some possibilities for choosing the function $g=g(x,y)$
in the Cauchy problem~\eqref{eq:02p} on concrete examples.

\begin{description}

\item{$1^\circ$}. The special case
$$
g=f
$$
is equivalent to the hodograph transformation with
an additional shift of the dependent variable, which gives $\xi=y-y_0$.

\item{$2^\circ$}. Setting
$$
g=\sqrt{1+f^2},
$$
we arrive at the method of the arc-length transformation~\cite{mor1979}.
In this case, the Cauchy problem~\eqref{eq:02p} takes the form
\begin{equation}
\begin{gathered}
x^{\prime}_\xi=\frac 1{\sqrt{1+f^2(x,y)}},\quad \ y^{\prime}_\xi=\frac {f(x,y)}{\sqrt{1+f^2(x,y)}};\\
x(0)=x_0,\quad \ y(0)=y_0.
\end{gathered}
\label{12.1.10.34*a}
\end{equation}

\item{$3^\circ$}. Choosing
$$
g=1+|f|,
$$
we obtain the Cauchy problem
\begin{equation}
\begin{gathered}
x^{\prime}_\xi=\frac 1{1+|f(x,y)|},\quad \ y^{\prime}_\xi=\frac {f(x,y)}{1+|f(x,y)|};\\
x(0)=x_0,\quad \ y(0)=y_0.
\end{gathered}
\label{12.1.10.34*c}
\end{equation}
Note that we use here the absolute value sign to generalize the results, since
the system~\eqref{12.1.10.34*c} can also be used in the case $f<0$ for numerical
integration of the problems having solutions with a root singularity~\cite{als2005}.

\item{$4^\circ$}. We can also take the function
$$
g=c_1+\bl(c_2+|f|^s\br)^{1/s}
$$
for $c_1\ge 0$, $c_2\ge 0$ ($|c_1|+|c_2|\not=0$), and $s>0$,  which is a generalization of the functions in
Items~$2^\circ$ and $3^\circ$.

\item{$5^\circ$}. A very convenient problem for analysis can be obtained if
we take
\begin{equation}
g=f/y
\label{eq:exp}
\end{equation}
in~\eqref{eq:02p}. In this case, the second equation of the system is immediately integrated
and, taking into account the initial condition, we get $y=y_0e^\xi$.
In addition, the variable~$x$ tends exponentially rapidly to a blow-up point~$x_*$ with increasing $\xi$. This
transformation will be called the special exp-type transformation.

\end{description}

\begin{remark}
From Items~$1^\circ$ and $2^\circ$ it follows that the method based on the hodograph
transformation and the method of the arc-length transformation are particular
cases of the method based on a non-local transformation of the general form~\eqref{*}.
\end{remark}

\begin{remark}
The functions $g$ in Items~$1^\circ\hbox{--}4^\circ$ correspond to the value
$k=1$ in~\eqref{eq:02q}, and the function~$g$ in Item~$5^\circ$ gives $k=\infty$.
\end{remark}

\begin{remark}
Nonlocal transformations of a special form were
used in \cite{mur2010,kudr2016a,kudr2016b} to obtain exact solutions and to linearize
some second-order ODEs.
\end{remark}

\subsection{Test problems and numerical solutions}\label{ss:3.2}

\begin{figure}
\centering
{\includegraphics[scale=0.34]{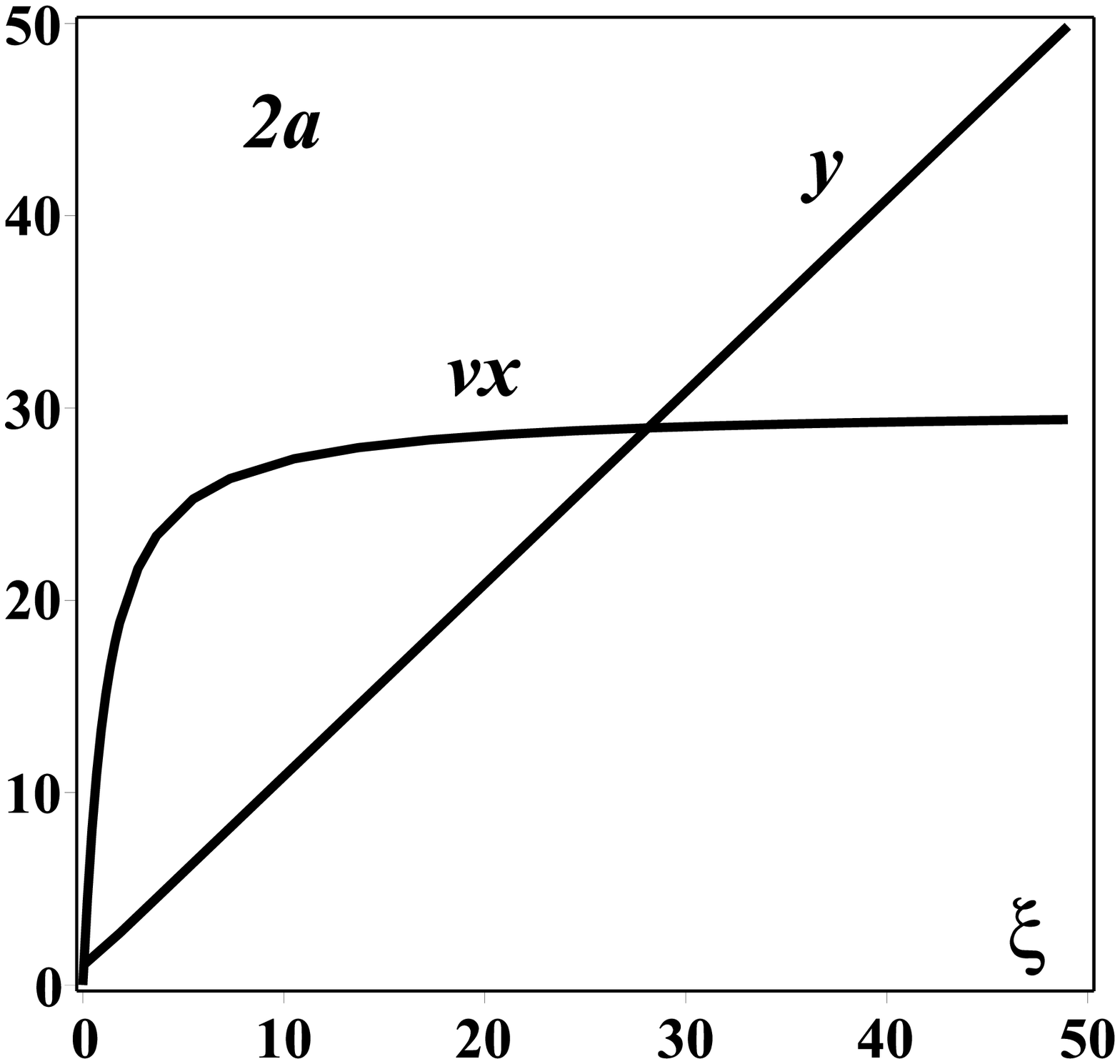}\ \includegraphics[scale=0.34]{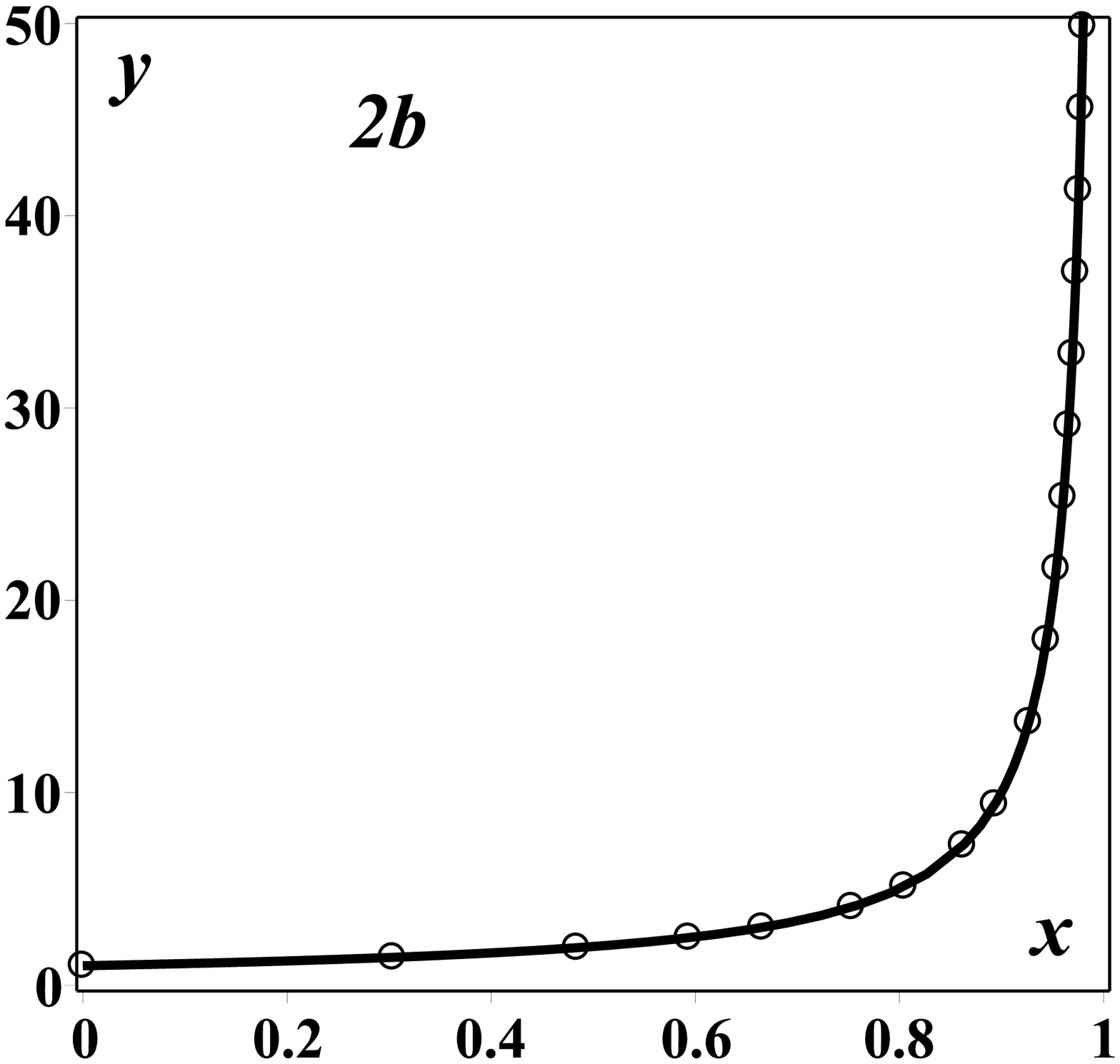}}
\caption{{\it 2a} \dash the dependences $x=x(\xi)$ and $y=y(\xi)$, obtained by numerical solution of the problem~\eqref{12.1.10.34*ab}
($\nu=30$); {\it 2b} \dash exact solution~\eqref{eq:02g} (solid line) and numerical solution
of the problem~\eqref{12.1.10.34*ab} (circles).}
\label{fig:Fig2}
\end{figure}

\begin{example}
For the test Cauchy problem \eqref{eq:02f} with $f=y^2$, the equivalent
problem for the system of equations~\eqref{12.1.10.34*a} takes the form
\begin{equation}
\begin{gathered}
x^{\prime}_\xi=\frac 1{\sqrt{1+y^4}},\quad \ y^{\prime}_\xi=\frac {y^2}{\sqrt{1+y^4}};\qquad
x(0)=0,\quad \ y(0)=1.
\end{gathered}
\label{12.1.10.34*ab}
\end{equation}
The second equation of this system is an equation with separable variables
whose solution is not expressed in elementary functions.

The numerical solution of the Cauchy problem~\eqref{12.1.10.34*ab}
in parametric form and its comparison with the exact solution~\eqref{eq:02g} are shown
in Fig.~\ref{fig:Fig2}.
\end{example}

\begin{example}
For the test Cauchy problem~\eqref{eq:02f}, the equivalent
the problem for the system of equations~\eqref{12.1.10.34*c} admits an exact solution,
which is expressed in terms of elementary functions in a parametric form as follows:
\begin{equation}
\begin{gathered}
x=1+\tfrac 12\xi-\tfrac12\sqrt{\xi^2+4},\quad \ y=\tfrac 12\xi+\tfrac12\sqrt{\xi^2+4}\qquad (\xi\ge 0).
\end{gathered}
\label{eq:CPExactPar}
\end{equation}
This solution satisfies the initial conditions $x(0)=0$ and $y(0)=1$ and has no singularities.
The function $x(\xi)$ is bounded, increases monotonically, and tends to its limiting value
$\ds x_*=\lim_{\xi\to\infty}x(\xi)=1$.
The function $y(\xi)$ increases monotonically and tends to infinity as $\xi\to\infty$.
At large $\xi$ we have $x\approx 1-\xi^{-1}$ and $y\approx \xi+\xi^{-1}$.

The curves $x=x(\xi)$ and $y=y(\xi)$, determined by the exact solution~\eqref{eq:CPExactPar}
(and also the curves obtained by numerical integration of the corresponding system~\eqref{12.1.10.34*c}
with $f=y^2$, $x_0=0$, and $y_0=1$), are very close to the curves shown in Fig.~\ref{fig:Fig2}
(they almost merge with them and therefore are not presented here).
\end{example}

\begin{example}
Consider the test problem~\eqref{eq:xx02}--\eqref{eq:xx03}, where $f=by^{\gamma}$, and take $g=f/y=by^{\gamma-1}$
(see Item~$5^\circ$ in Section~3.1).
Substituting these functions into~\eqref{eq:02p}, we obtain the Cauchy problem
\begin{equation}
\begin{gathered}
x^{\prime}_\xi=\frac 1{by^{\gamma-1}},\quad \ y^{\prime}_\xi=y\quad \ (\xi>0);\\
x(0)=0,\quad y(0)=a,
\end{gathered}
\label{eq:02r}
\end{equation}
where $a>0$, $b>0$, and $\gamma>1$.
The exact solution of the problem \eqref{eq:02r} is written as follows:
\begin{equation}
x=\frac 1{a^{\gamma-1}b(\gamma-1)}\bl[1-e^{-(\gamma-1)\xi}\br],\quad \ y=ae^\xi.
\label{eq:02pe}
\end{equation}
It can be seen that the unknown function $x=x(\xi)$ tends exponentially to the asymptotic value
$\ds x_*=\frac 1{a^{\gamma-1}b(\gamma-1)}$ as $\xi\to \infty$.

Let $b=1$ and $\gamma=2$. The numerical solutions of the problems~\eqref{eq:02h} and~\eqref{eq:02r},
obtained by the Runge--Kutta method of the fourth-order of approximation,
are shown in Fig.~\ref{fig:Fig3} for $a=1$ and $a=2$ and the same step of integration, equal to~$0.2$.
We note that for this stepsize, the maximum difference between the exact solution~\eqref{eq:02pe}
and the numerical solution of the system~\eqref{eq:02r} is $0.0045\%$
(and for stepsize $0.4$, respectively, $0.061\%$).

\medskip
\begin{figure}
\centering
{\includegraphics[scale=0.34]{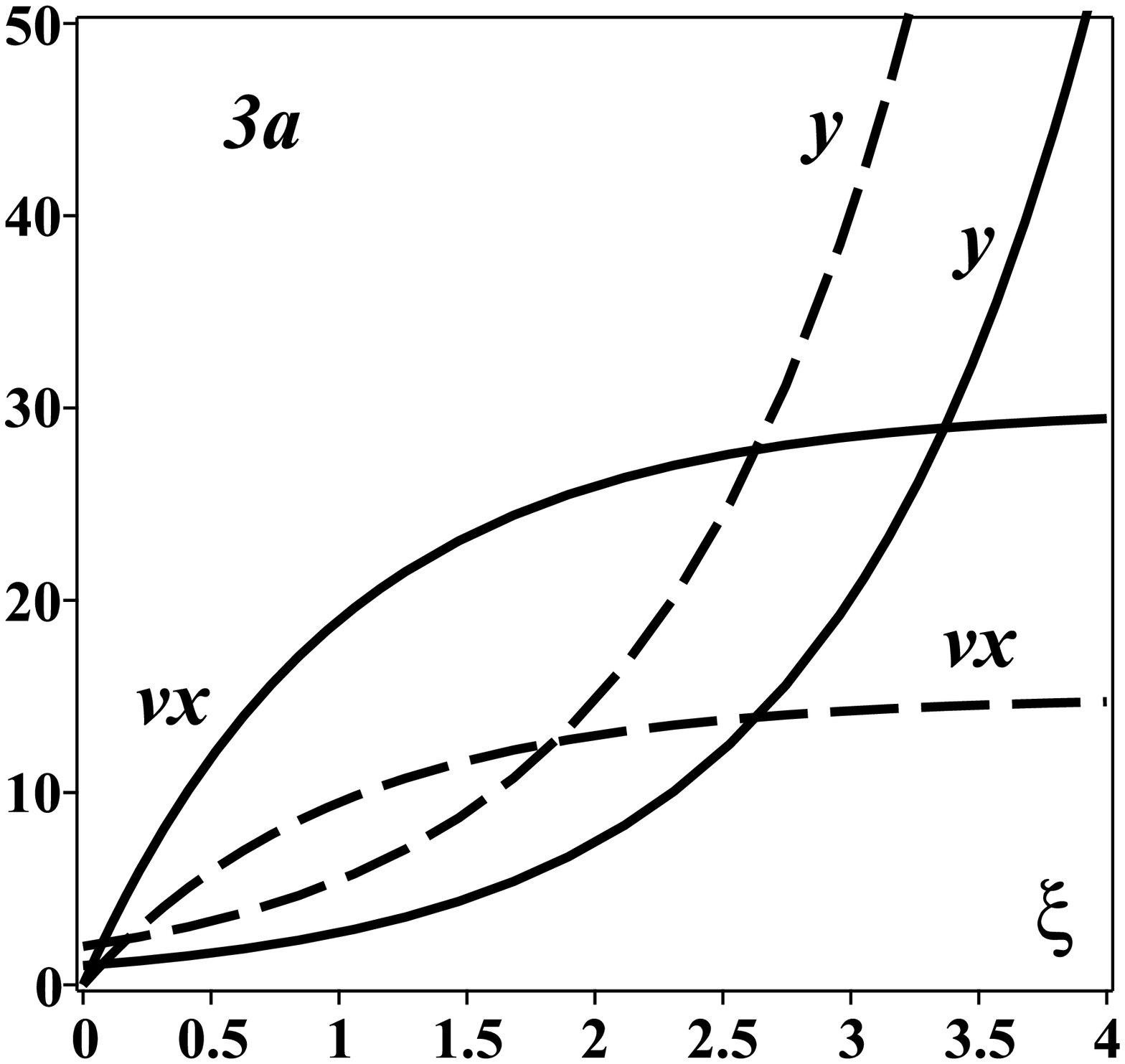}\ \includegraphics[scale=0.34]{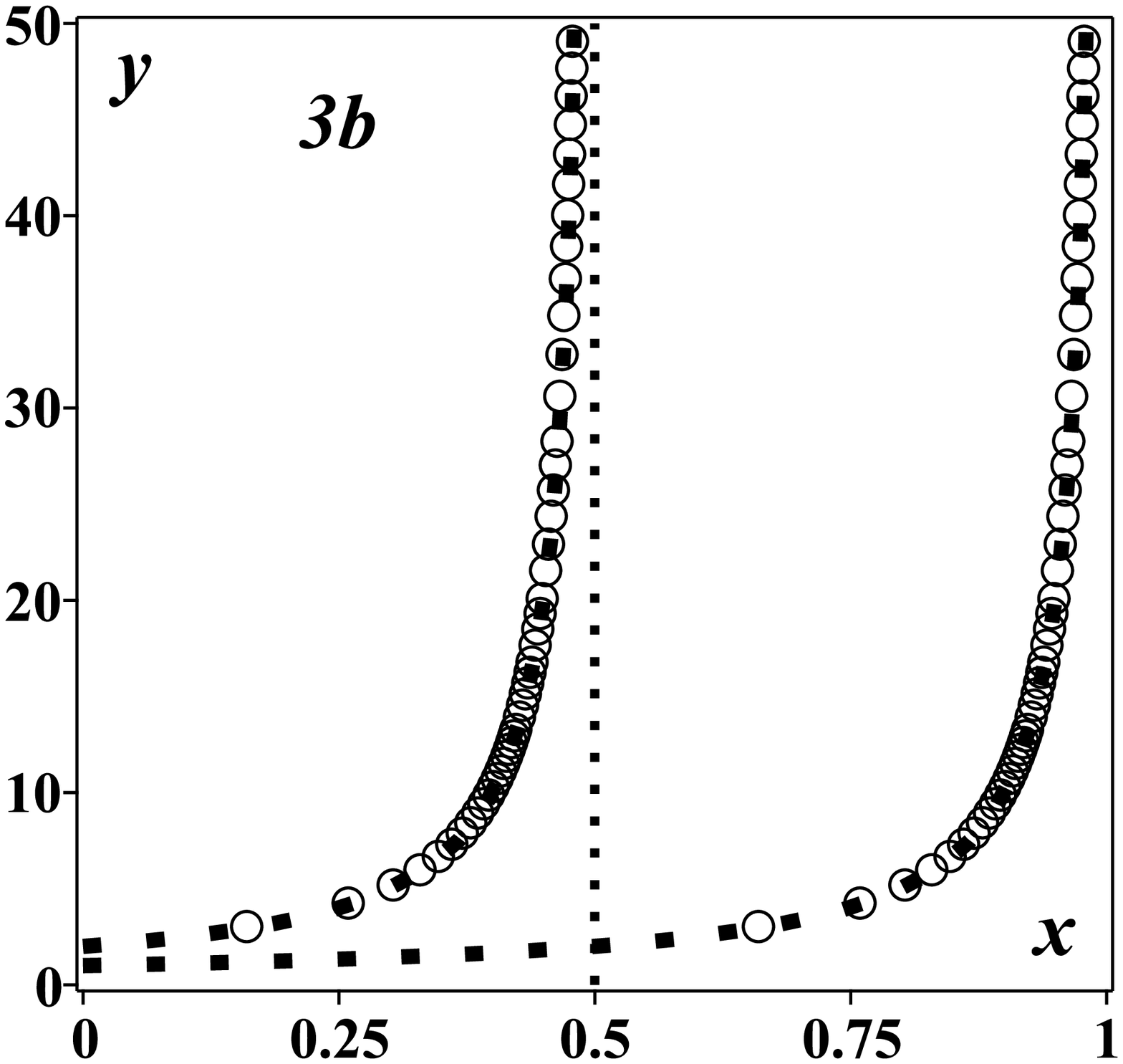}}
\caption{{\it 3a} \dash the dependences $x=x(\xi)$ and $y=y(\xi)$, obtained by numerical solution
 of the problem~\eqref{eq:02r} for $b=1$, $\gamma=2$ with $a=1$ (solid lines) and $a=2$ (dashed lines) ($\nu=30$);
{\it 3b} \dash numerical solution of the problem~\eqref{eq:02h} for $b=1$, $\gamma=2$ (circles)
and numerical solution of the problem~\eqref{eq:02r} for $b=1$, $\gamma=2$ (points);
for left curves $a=2$ and for right curves $a=1$.}
\label{fig:Fig3}
\end{figure}

It can be seen (see Figs.~\ref{fig:Fig1} and \ref{fig:Fig3}) that the numerical solutions
are in a good agreement, but the rates of their approximation to the required asymptote $x=x_*$
are significantly different. For example, for the system~\eqref{eq:02h},
in order to obtain a good approximation to the asymptote,
it is required to consider the interval $t\in[1,2980]$, and for the system~\eqref{eq:02r} it suffices to take
$\xi\in[0,4]$. Therefore, there is reason to believe that the method described in Item~$5^\circ$
(a special case of the transformation~\eqref{*}) is more
efficient than the method based on the differential transformation
(see Section~2.1).

For comparison, similar calculations were also performed applying
the method based on the hodograph transformation (see Section~3.1, Item~$1^\circ$),
and the method of the arc-length transformation (see Section~3.1, Item~$2^\circ$).
For both of these methods, in order to obtain a good approximation to the asymptote,
it is required to consider the interval $\xi\in[0,49]$.
To control a numerical integration process, the calculations were carried out
with the aid of the three most important and powerful
mathematical software packages: Maple (2016), Mathematica~(11), and MATLAB (2016a).
It was found that the method based on the use of a special case of the system~\eqref{eq:02p}
with $g=f/y$ (see Item~$5^\circ$) is essentially
more efficient than the method based on the hodograph transformation
and the method of the arc-length transformation.
\end{example}
\medskip


\begin{example}
We now consider the Cauchy problem~\eqref{*12.1.10.33} for $b=1$, which
is determined by the exponential $f=e^y$. Substituting  the function
$g=f/y=e^y/y$ (see Item~$5^\circ$ in Section~3.1) into the system~\eqref{eq:02p}, we obtain
\begin{equation}
\begin{gathered}
x^{\prime}_t=ye^{-y},\quad \ y^{\prime}_t=y;\\
x(0)=0,\quad y(0)=a.
\end{gathered}
\notag
\end{equation}
The exact solution of this problem in a parametric form is defined by the formulas
\begin{align}
x=e^{-a}-\exp(-a\,e^t),\quad \ y=ae^t\qquad (t\ge 0),
\notag
\end{align}
which do not have singularities. The function $x=x(t)$ is bounded,
monotonically increases with increasing~$t$ and very rapidly tends to the asymptote
$\ds x_*=\lim_{t\to\infty}x(t)=e^{-a}$,
and the function $y=y(t)$ is unbounded and grows exponentially with respect to~$t$.
\end{example}

\subsection{Generalizations based on the use of differential constraints}\label{ss:3.3}

Let us show that the method based on introducing a non-local variable~\eqref{*}
allows a further generalization.

We add to the equation~\eqref{eq:02} a first-order differential constraint~\cite{pol2012} of the form
\begin{align}
\xi'_x=g(x,y,\xi)
\label{eq:1000*}
\end{align}
and the initial condition $\xi(x=x_0)=\xi_0$.

In a particular case, when the function~$g$ does not depend on~$\xi$,
the use of the differential constraint~\eqref{eq:1000*}, after integrating it over~$x$,
leads to the non-local variable~\eqref{*} for $\xi_0=0$
(therefore, the method based on the differential constraint generalizes the method based on
introducing a non-local variable).

From \eqref{eq:02} and \eqref{eq:1000*} we obtain the following
system of ordinary differential equations:
\begin{align}
x^{\prime}_\xi=\frac 1{g(x,y,\xi)},\quad \ y^{\prime}_\xi=\frac {f(x,y)}{g(x,y,\xi)}.
\label{eq:1001}
\end{align}

In a particular case,
$$
g(x,y,\xi)=f_x+\xi f_y,
$$
the system \eqref{eq:1001} coincides with the system~\eqref{eq:02dd}, in which the variable~$t$
must be redenoted by~$\xi$.
If, in addition, we set $\xi_0=f(x_0,y_0)$, then, up to renaming of variables,
we also obtain the initial conditions~\eqref{eq:02e}.
It follows that the method based on the differential constraint of
general form~\eqref{eq:1000*} generalizes the method based on introducing a differential variable (see Section~2.1).

\subsection{Comparison of efficiency of various transformations for numerical integration
of first-order ODE blow-up problems}\label{ss:3.4}

In Table~1, a comparison of the efficiency of the numerical integration methods,
based on various non-local transformations of the form~\eqref{*} and differential
constraints of the form~\eqref{eq:1000*} is presented by using the example of the test blow-up problem
for the first-order ODE~\eqref{eq:02f} with $f=y^2$.
The comparison is based on the number of grid points needed to perform calculations with the same
maximum error (approximately equal to $0.1$,~$0.01$, and $0.005$).
In the last line of Table~1 for Example~4 we take $a=b=1$ and $s=2$.



\begin{table}[!ht]
\begin{center}
\begin{tabular}{|l|l|r|r|r|l|l|r|r|r|}
\hline
\multicolumn{5}{|c|}{\vphantom{$\frac12$} {\small Error$_{{\rm max}}, \% =0.1$} }      \\ \hline \hline
{\small Transformation or}     &  \quad {\small Function $g$}  &  {\small Max. interval}  & {\small Stepsize}      & {\small Grid points} \\
{\small differential constraint}  &  \quad {\small or Example}  &  {\small $\xi_{{\rm max}}\qquad$}  & {\small $h\quad$}  & {\small number}  $N$  \\  \hline
{\small Hodograph, Item~$1^\circ$}         &  $g{=}f$                                      &  48.99 &   0.2300          & 213 \\ \hline
{\small Arc-length, Item~$2^\circ$}        &  {\small $g{=}\sqrt{1{+}f^2}$}                &  49.20 &   0.3000          & 164 \\ \hline
{\small Nonlocal, Item~$3^\circ$}  &  {\small $g{=}1{+}|f|$}                               &  50.00 &   0.4000          & 125 \\ \hline
{\small Special exp-type, Item~$5^\circ$} &  $g{=}f/y$                                     &  3.925 &   0.1570          & 25  \\ \hline
{\small Diff. constraint}  &  $g{=}f/[y(1+2\xi)]$                                          &  1.543 &   0.0643          & 24  \\ \hline
{\small Differential, modified}        &  {\small Example 4 with $\lambda=2$}              &  3.910 &   0.2300          & 17  \\ \hline
\hline
\multicolumn{5}{|c|}{\vphantom{} {\small Error$_{{\rm max}}, \% =0.01$} }     \\ \hline \hline
{\small Transformation or}     &  \quad {\small Function $g$}  &  {\small Max. interval}  & {\small Stepsize}      & {\small Grid points} \\
{\small differential constraint}  &  \quad {\small or Example}  &  {\small $\xi_{{\rm max}}\qquad$}  & {\small $h\quad$}  & {\small number}  $N$  \\  \hline
{\small Hodograph, Item~$1^\circ$}         &  $g{=}f$                                      & 49.01  &     0.130          &  377  \\ \hline
{\small Arc-length, Item~$2^\circ$}        &  {\small $g{=}\sqrt{1{+}f^2}$}                & 49.30  &     0.170          &  290  \\ \hline
{\small Nonlocal, Item~$3^\circ$}  &  {\small $g{=}1{+}|f|$}                               & 50.14  &     0.230          &  218  \\ \hline
{\small Special exp-type, Item~$5^\circ$} &  $g{=}f/y$                                     & 3.960  &     0.090          &  44   \\ \hline
{\small Diff. constraint}  &  $g{=}f/[y(1+2\xi)]$                                          & 1.540  &     0.035          &  44   \\ \hline
{\small Differential, modified}        &  {\small Example 4 with $\lambda=2$}              & 3.900  &     0.130          &  30   \\ \hline
\hline
\multicolumn{5}{|c|}{\vphantom{} {\small Error$_{{\rm max}}, \% =0.005$} }     \\ \hline \hline
{\small Transformation or}     &  \quad {\small Function $g$}  &  {\small Max. interval}  & {\small Stepsize}      & {\small Grid points} \\
{\small differential constraint}  &  \quad {\small or Example}  &  {\small $\xi_{{\rm max}}\qquad$}  & {\small $h\quad$}  & {\small number}  $N$  \\  \hline
{\small Hodograph, Item~$1^\circ$}         &  $g{=}f$                                      & 49.035  &     0.1050         &  467  \\ \hline
{\small Arc-length, Item~$2^\circ$}        &  {\small $g{=}\sqrt{1{+}f^2}$}                & 49.266  &     0.1380         &  357  \\ \hline
{\small Nonlocal, Item~$3^\circ$}  &  {\small $g{=}1{+}|f|$}                               & 50.135  &     0.1850         &  271  \\ \hline
{\small Special exp-type, Item~$5^\circ$} &  $g{=}f/y$                                     & 3.9150  &     0.0725         &  54   \\ \hline
{\small Diff. constraint}  &  $g{=}f/[y(1+2\xi)]$                                          & 1.5420  &     0.0291         &  53   \\ \hline
{\small Differential, modified}        &  {\small Example 4 with $\lambda=2$}              & 3.9140  &     0.1030         &  38   \\ \hline
\end{tabular}
\caption{Various types of analytical transformations applied for numerical
integration of the problem~\eqref{eq:02f} for $f=y^2$
with a given accuracy (percent errors are $0.1$, $0.01$, and $0.005$ for
$\Lambda_\text{m}\le 50$) and their basic parameters  (maximum interval,
stepsize, grid points number).}
\label{tab:Transfs-ODE1}
\end{center}
\end{table}

It can be seen that
for the first three transformations it is necessary to use a lot of grid points (the hodograph transformation is the least effective).
This is due to the fact that in these cases $x$~tends to the point~$x_ *$ rather slowly for large~$\xi$ ($x_*-x\sim 1/\xi$, $y\sim \xi$).
The last three transformations require a significantly less number of grid points; in these cases $x$~tends exponentially fast
to the point~$x_ *$ for large~$\xi$.
In particular, the use of the exp-type transformation with $g=t/y$ gives rather good results.
The most effective analytical transformation is a modification of the method of differential transformations (see Example 4 in Section~2.3).

\subsection{Complex blow-up problems, in which the right-hand side of the equation and the solution can change the sign}\label{ss:3.5}

Up to now it has been assumed that the right-hand side of the equation~\eqref{eq:02} is positive.
If we remove this restriction and assume that the function $f=f(x,y)$, and also the solution $y=y(x)$,
can change sign for $x>x_0$ (in addition, $y_0$ can be any sign),
then we can act in the following two stages:

\begin{description}

\item{(i)} In the first stage, the problem \eqref{eq:02}--\eqref{eq:02a} is integrated by applying
a standard numerical method and the value $\Lambda_m=\min\{|y/y_0|,\ |y^{\prime}_x/y|\}$ (for $y_0\sim 1$) is calculated at the grid points~$x_m$.
When the value $\Lambda_\text{m}$ becomes sufficiently large (for example, $\Lambda_\text{m}=30$),
then the calculation stops. The results of the integration performed are used in the region $x_0\le x\le x_m$.

\item{(ii)} In the second stage, we set the initial condition $y(x_m)=y_m$ at the point~$x_m$
to obtain a solution in the region $x>x_m$ (here $y_m$ is the value calculated on the first stage).
Then a non-local transformation is applied (for example, described in Item~$5^\circ$ in Section~3.1)
and the resulting problem is integrated applying standard fixed-step numerical methods.

\end{description}

\section{Problems for first-order equations, the right-hand side of which has\\ singularities or zeros}\label{ss:4}

\subsection{Blow-up problems for equations, the right-hand side of which has\\ singularities in~$x$}\label{ss:4.1}

In this section we will analyze several blow-up problems for equations of the
form~\eqref{eq:02}, the right-hand side of which has a singularity at some
$x=x_\text{s}$, i.e. $\ds\lim_{x\to x_\text{s}}f(x,y)=\infty$.

We assume that the right-hand side of equation~\eqref{eq:02} can be
represented as a product of two functions
\begin{equation}
f(x,y)=f_\text{b}(x,y)f_\text{s}(x,y), \label{eeqq:03}
\end{equation}
where the function $f_\text{b}$ has the same properties as the function~$f$ in
Section~2.1 (i.e. that the problem \eqref{eq:02}--\eqref{eq:02a},
where the function~$f$ is replaced by $f_\text{b}$, has a blow-up solution).

Moreover, we will assume that the function $f_\text{s}$ has an integrable or non-inte\-grable singularity at
$x=x_\text{s}$, so that $\ds\lim_{x\to x_\text{s}}f(x,y)=\infty$, and $f_\text{s}>0$ at $x_0<x_\text{s}$.

It is interesting to see how the two singularities of this problem will interact:
the blow-up singularity and the coordinate singularity at $x=x_\text{s}$.

For the sake of clarity, we give the following test problems and illustrative examples.

\begin{example}
Consider the two-parameter test Cauchy problem:
\begin{equation}
y^{\prime}_x=\frac{y^2}{b-x}; \quad \ y(0)=a, \label{eeqq:03a}
\end{equation}
where $a>0$ and $b>0$.
For this problem, we have $f_\text{b}=y^2$ and
$f_\text{s}=1/(b-x)$. The right-hand side of equation~\eqref{eeqq:03a} has a
pole of the first order at the point $x=x_\text{s}=b$ (i.e., there exists a
non-integrable singularity at this point); and the right-hand side of the
equation becomes negative if $x>b$.

The exact solution of the problem~\eqref{eeqq:03a} has the form
\begin{equation}
y=\BL[\ln\Bl(1-\frac xb\Br)+\frac 1a\BR]^{-1}. \label{eeqq:03b}
\end{equation}

The singular point of this solution is determined by the formula
\begin{equation}
x_*=b\bl(1-e^{-1/a}\br)<b. \label{eeqq:03c}
\end{equation}
Here the blow-up singularity ``overtakes" the non-integrable singularity of the equation
at the point $x_\text{s}=b$. If $a\to\infty$, we have $x_*\to b$.
\end{example}

\begin{example}
Consider the test Cauchy problem with a stronger coordinate singularity:
\begin{equation}
y^{\prime}_x=\frac{y^2}{(b-x)^2}; \quad \ y(0)=a, \label{eeqq:03d}
\end{equation}
where $a>0$ and $b>0$. For this problem, we have $f_\text{b}=y^2$ and
$f_\text{s}=1/(b-x)^2$. The right-hand side of equation~\eqref{eeqq:03d} has
a pole of the second order at the point $x_\text{s}=b$ (i.e., there exists
a non-integrable singularity at this point); and right-hand side of this equation is positive for all~$x$.

The exact solution of the problem~\eqref{eeqq:03d} is defined by the formula
\begin{equation}
y=\frac {ab}{a+b}\BL[1+\frac {ab}{b^2-(a+b)x}\BR]. \label{eeqq:03e}
\end{equation}

The blow-up point is determined as
\begin{equation}
x_*=\frac {b^2}{a+b}<b. \label{eeqq:03f}
\end{equation}
Here, as in Example~9, the blow-up singularity ``overtakes" the non-integrable
singularity of the equation at the point $x_\text{s}=b$.
\end{example}

\begin{remark}
A qualitatively similar picture will occur also for the
problem \eqref{eq:02}--\eqref{eq:02a} with $f=y^2/\sqrt{b-x}$, which has
an integrable singularity at the point $x=b$. In this problem, the right-hand
side of the equation is defined only on a part of the $x$-axis.
\end{remark}

The solution property, described in Examples~9 and~10, has a
general characteristic. Namely, let us assume that the right-hand side of
equation~\eqref{eq:02} has the form~\eqref{eeqq:03}, where the functions
$f_\text{b}$ and $f_\text{s}$ have the properties described at the beginning
of this section. Then the problem \eqref{eq:02}--\eqref{eq:02a} has a
blow-up solution, and the domain of definition of this solution is located to
the left of the point $x_\text{s}$ (i.e., $x_*<x_\text{s}$, where $x_*$ is
the blow-up point).

The methods described in Sections~2 and 3 can be applied
for solving this type of problems with the coordinate singularity.

\begin{figure}
\centering
\vskip-2pc
{\includegraphics[scale=0.34]{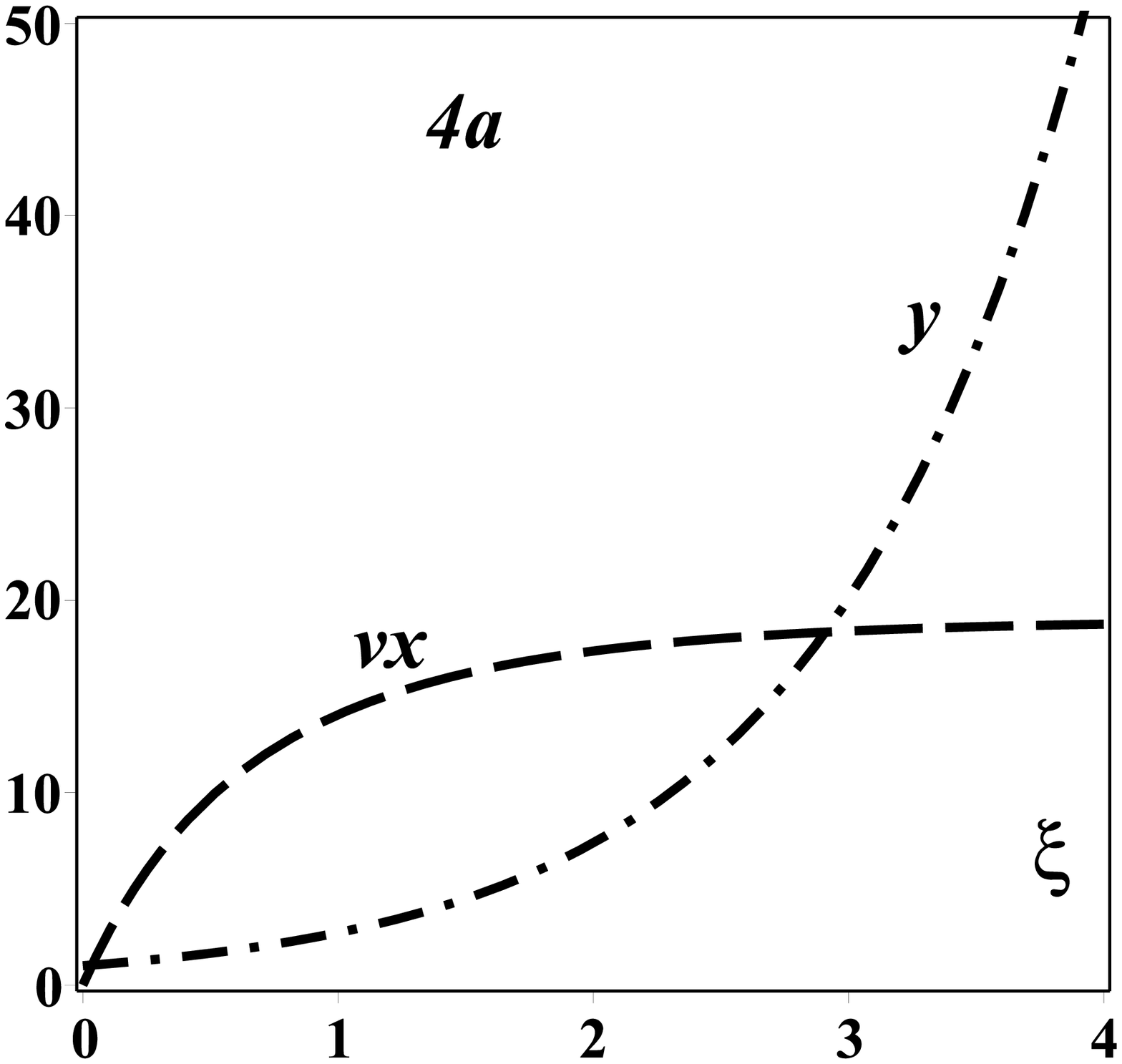}\ \includegraphics[scale=0.34]{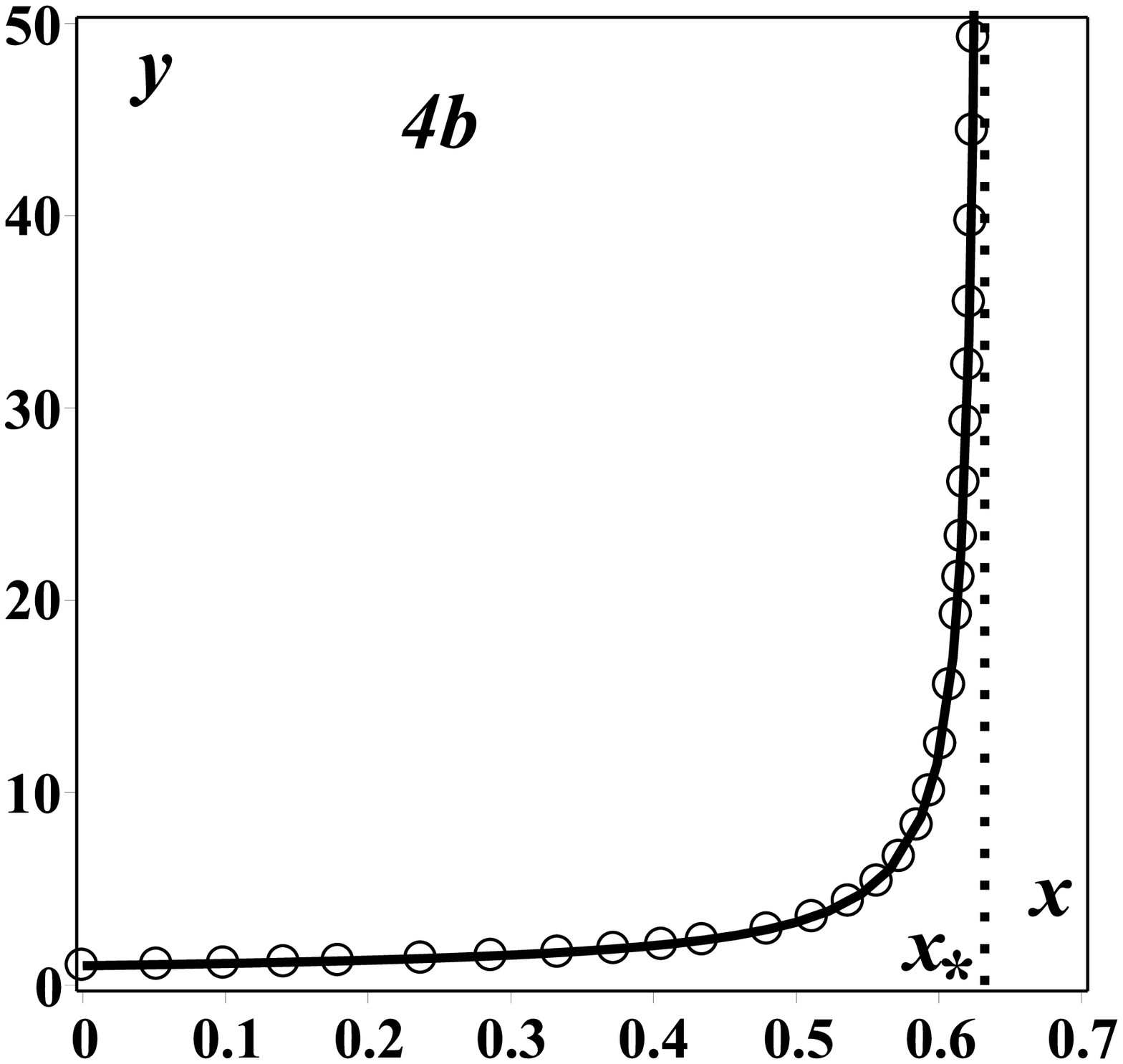}}
\vskip-1pc
\caption{{\it 4a}\dash numerical solutions of system~\eqref{eeqq:03g}, the dependences $x=x(\xi)$ and $y=y(\xi)$
(where $\nu=30$); {\it 4b}\dash exact solution~\eqref{eeqq:03b}, solid line, and the numerical solution of
system~\eqref{eeqq:03g} ($a=b=1$ and $x_*=0.6321$).}
\label{fig:Fig7}
\end{figure}

\begin{example}
For test problem \eqref{eeqq:03a} with $f=y^2/(b-x)$,
we take (see Item~$5^\circ$ in~Section~3.1),
$$
g=\frac{f}{y}=\frac{y}{b-x}.
$$
By substituting these functions in~\eqref{eq:02p}, we arrive at the Cauchy
problem for a system of coupled ODEs
\begin{equation}
\begin{gathered}
x^{\prime}_\xi=\frac {b-x}{y},\quad \ y^{\prime}_\xi=y\quad \ \ (\xi>0);\\
x(0)=0,\quad y(0)=a.
\end{gathered}
\label{eeqq:03g}
\end{equation}
The exact solution of this problem reads
\begin{equation}
x=b\BL\{1-\exp\BL[-\frac 1a(1-e^{-\xi})\BR]\BR\},\quad \ y=ae^\xi.
\label{eeqq:03h}
\end{equation}

The numerical solution of the problem~\eqref{eeqq:03g} with $a=b=1$ is
presented in Fig.~\ref{fig:Fig7}{\sl a}; the dependences $x=x(\xi)$ and $y=y(\xi)$
are obtained by the fourth-order Runge--Kutta method.
Fig.~\ref{fig:Fig7}{\sl b} shows a comparison of the exact solution~\eqref{eeqq:03b} of
the Cauchy problem~\eqref{eeqq:03a} for one equation with the numerical
solution of the problem for the system of two equations~\eqref{eeqq:03g}.
\end{example}

\begin{remark}
A form of the right-hand side of equation~\eqref{eq:02}
with a coordinate singularity at the point $x_\text{s}$ can mislead the
researcher, inexperienced in blow-up problems. As a result, the researcher will start to
refine a mesh (making it thinner) in the neighborhood of the point
$x_\text{s}$ (what should not to do).
\end{remark}

\subsection{Blow-up problems for equations, the right-hand side of which has zeros}\label{ss:4.2}

In this section we will analyze blow-up problems for equations of the
form \eqref{eq:02}, the right-hand side of which vanishes at some
$x=x_\text{z}$, i.e. $f(x_\text{z},y)=0$.

Let us assume that the right-hand side of equation~\eqref{eq:02} can be
represented as a product of two functions
\begin{equation}
f(x,y)=f_\text{b}(x,y)f_\text{z}(x,y), \label{eeqq:03i}
\end{equation}
where the function~$f_\text{b}$  has the same properties as the function~$f$ in
Section~2.1 (i.e., the problem \eqref{eq:02}--\eqref{eq:02a},
where the function~$f$ is replaced by~$f_\text{b}$, has a blow-up solution).
Moreover, we will assume that the function~$f_\text{z}$ vanishes at
$x=x_\text{z}$, so that $f(x_\text{z},y)=0$, and $f_\text{z}>0$ at
$x_0<x_\text{z}$.

\begin{example}
Consider the test two-parameter Cauchy problem
\begin{equation}
y^{\prime}_x=y^2(b-x); \quad \ y(0)=a, \label{eeqq:03j}
\end{equation}
where $a>0$ and $b>0$.
For this problem, we have $f_\text{b}=y^2$ and
$f_\text{z}=b-x$. The right-hand side of equation~\eqref{eeqq:03j} becomes
zero at the point $x=x_\text{z}=b$; and the right-hand side of the equation
becomes negative if $x>b$.

It is interesting to see how two features of different types of such problem
will interact: on the one hand, a possible blow-up singularity
(which leads to an unlimited growth of the right-hand side of the equation), and on the other hand,
vanishing  of the right-hand side of the equation at $x=x_\text{z}$.

The exact solution of the problem~\eqref{eeqq:03j} has the form
\begin{equation}
y=\frac a{\frac12ax^2-abx+1}. \label{eeqq:03k}
\end{equation}
The existence or absence of a blow-up singularity in this solution is determined
by the existence or absence of real roots of the quadratic equation $\frac12ax^2{-}abx{+}1{=}0$.

\medskip
The elemental analysis shows that there are two qualitatively different cases:

\begin{figure}
\centering
\vskip-2pc
{\includegraphics[scale=0.4]{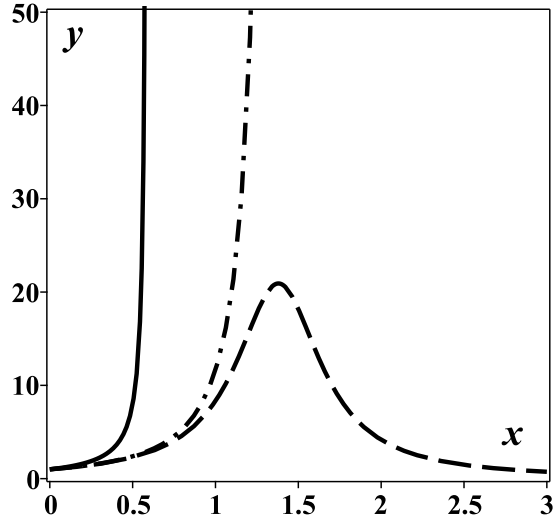}}
\vskip-1pc
\caption{Exact solutions~\eqref{eeqq:03k} of the Cauchy
problem~\eqref{eeqq:03j} for various values of the parameters: $a=1$, $b=1.38$
(dashed line); $a=1$, $b=2$ (solid line); and $a=1$, $b=\sqrt{2}$ (dashed-dot line).}
\label{fig:Fig8}
\end{figure}

\begin{description}

\item{(i)} If $0<b<\sqrt{2/a}$, there exists a smooth continuous solution of the problem for all $x\ge0$.
It  is monotonically increasing on the interval $0\le x<b$, reaches the maximum
value $\ds y_m=\frac a{1-\frac12 ab^2}$, and decreases for $x>b$.

\item{(ii)} If $b\ge \sqrt{2/a}$, formula~\eqref{eeqq:03k} defines a monotonically
increasing blow-up solution with the singular point
\begin{equation}
\ds x_*=b-\sqrt{b^2-\frac 2a}. \label{eeqq:03l}
\end{equation}

\end{description}

Therefore in this example the existence or absence of a blow-up solution is
determined by a simple relation between the parameters $a$ and~$b$: if $b\ge
\sqrt{2/a}$, then there exists a blow-up solution, otherwise there is no a
blow-up solution. The value $b=\sqrt{2/a}$ is a point of bifurcation of
the two-parameter problem~\eqref{eeqq:03j}.

The exact solutions of the problem~\eqref{eeqq:03j} obtained by
formula~\eqref{eeqq:03k} are presented in Fig.~\ref{fig:Fig8}\ for various values of
the parameters: $a=1$ and $b=1$, $b=2$, and $b=\sqrt{2}$ (the critical value
at which there exists blow-up solutions).

Similar problems are complicated for numerical integration in a wide range of
changing of the parameters $a$ and~$b$. First, it is reasonable to produce a
direct integration of the problem~\eqref{eeqq:03j} by the Runge--Kutta
methods for specific values of $a$ and~$b$. If a rapid growth of the
solution occurs, then  the methods described in
Sections~2 or 3 should be applied (see also Section~3.4).

More efficient methods allowing to integrate numerically similar and
more complex non-monotonic blow-up problems
will be described in subsequent publications.

\end{example}

\begin{remark}
A formal replacing $b-x$ to $(b-x)^2$ in the
problem~\eqref{eeqq:03j} leads to a blow-up solution for any $a>0$ and $b>0$.
\end{remark}

\section{Problems for first-order equations. Two-sided estimates of\\ the critical value}\label{ss:5}

\subsection{Autonomous equations. Analytical formula for the critical value}\label{ss:5.1}

We consider the Cauchy problem for an autonomous equation of the general form
\begin{align}
y^{\prime}_x=f(y)\quad (x>0),\quad \ y(0)=a.
\label{12.1.1.1a}
\end{align}
We assume that $a>0$ and $f(y)>0$  is a continuous function that is defined
for all $y\ge a$. An exact solution of the Cauchy problem~\eqref{12.1.1.1a}
for $x>0$ can be represented implicitly as follows:
\begin{align}
x=\int^y_a\frac{d\xi}{f(\xi)}.
\label{*12.1.10.32a}
\end{align}
This solution is a blow-up solution if and only if
there exists a finite definite integral in~\eqref{*12.1.10.32a} for $y=\infty$.
In this case, the critical value~$x_*$ is calculated as follows:
\begin{align}
x_*=\int^\infty_a\frac{d\xi}{f(\xi)}.
\label{*12.1.10.32b}
\end{align}

Let the conditions formulated after the problem~\eqref{12.1.1.1a} be satisfied.

\textsl{A necessary criterion for the existence of a blow-up solution} is:
$$
\ds\lim_{y\to\infty} \frac{f(y)}{y}=\infty.
$$

\textsl{Sufficient criterion of the existence of a blow-up  solution}.
Let the conditions formulated above are also satisfied and the limiting ratio,
\begin{align}
\lim_{y\to\infty}\frac{f(y)}{y^{1+\kappa}}= s,\quad \ 0<s\le\infty,
\label{*12.1.10.32b**}
\end{align}
takes place for some parameter $\kappa>0$. Then the solution of the Cauchy problem~\eqref{12.1.1.1a} is a blow-up solution.

\medskip
If $f(y)$ is a differentiable function, then instead of~\eqref{*12.1.10.32b**}
we can propose an equivalent sufficient criterion of the existence of a blow-up  solution:
$$
\lim_{y\to\infty}\frac{f'_y(y)}{y^{\kappa}}= s_1,\quad \ 0<s_1\le\infty\qquad (\kappa>0).
$$

\subsection{Non-autonomous equations. One-sided estimates}\label{ss:5.2}

We consider the Cauchy problem for a first-order non-autonomous equation of the general form
\begin{align}
y^{\prime}_x=f(x,y)\quad (x>0),\quad \ y(0)=a.
\label{eq:*1}
\end{align}

We assume that $f(x,y)$ is a continuous function and the conditions
\begin{align}
f(x,y)\ge g(y)>0\quad \text{for all}\quad y\ge a>0, \quad \ \ x\ge 0
\label{*12.1.10.32d}
\end{align}
are satisfied. We also assume that there exists a finite integral
\begin{align}
I_g=\int^\infty_{a}\frac{d\xi}{g(\xi)}<\infty.
\label{*12.1.10.32e}
\end{align}
Then the solution $y=y(x)$ of the Cauchy problem~\eqref{eq:*1} is a blow-up solution,
and the critical value~$x_*$ satisfies the inequality
\begin{align}
x_*\le I_g.
\label{*12.1.10.32f}
\end{align}

This estimate follows from the inequality (see, for example, the corresponding
comparison theorems in~\cite{mikh1965,kam}):
\begin{align}
y(x)\ge y_g(x),
\label{*12.1.10.32f1}
\end{align}
where $y(x)$ is the solution of the Cauchy problem~\eqref{eq:*1}, and $y_g(x)$~is the solution
of the auxiliary Cauchy problem
\begin{align}
y^{\prime}_x=g(y)\quad \ (x>0),\quad \ y(0)=a.
\label{*12.1.10.32f2}
\end{align}

\begin{example}
We consider the Cauchy problem for the Abel equation of the first kind
\begin{align}
y^{\prime}_x=y^3+h(x)\quad (x>0);\quad \ y(0)=1.
\label{*12.1.10.32f*}
\end{align}
If $h(x)\ge 0$ for $x\ge 0$, then the inequality is valid
$$
f(x,y)\equiv y^3+h(x)\ge g(y)\equiv y^3>0\quad \text{for all}\quad y>1.
$$
Calculating the integral~\eqref{*12.1.10.32e} with $g(y)=y^3$, we obtain
\begin{align}
I_g=\int^\infty_{1}\frac{d\xi}{\xi^3}=\frac 1{2}<\infty.
\label{*12.1.10.32e*}
\end{align}
Therefore, the solution of the Cauchy problem~\eqref{*12.1.10.32f*} for $h(x)\ge 0$ is a blow-up solution,
and $x_*\le \frac12$.
\end{example}

\subsection{Non-autonomous equations. Two-sided estimates}\label{ss:5.3}

We consider two cases in which the one-sided estimate~\eqref{*12.1.10.32f} can be improved.
We introduce the notations
\begin{align}
I_1=\int^\infty_{y_0}\frac{d\xi}{f(0,\xi)},\quad \
I_2=\int^\infty_{y_0}\frac{d\xi}{f(I_1,\xi)}.
\label{*12.1.10.32hz}
\end{align}

\textsl{Case $1^\circ$.}  Let $f_x\ge 0$.
Suppose that the integral~$I_1$ in~\eqref{*12.1.10.32hz} exists and is finite.
Suppose also that the conditions,
\begin{align}
f(x,y)>0, \quad f_x(x,y)\ge 0 \quad \text{for all} \quad 0\le x\le I_1, \ \ y\ge y_0>0,
\label{*12.1.10.32g}
\end{align}
are satisfied. Then the integral~$I_2$ exists and the inequalities are valid:
\begin{equation}
f(0,y)\le f(x,y)\le  f(I_1,y) \quad \text{for}\quad 0\le x\le I_1
\label{*12.1.10.32h}
\end{equation}
and
\begin{align}
y_1(x)\le y(x)\le y_2(x) \quad \text{for}\quad 0\le x\le I_2\le I_1.
\label{*12.1.10.32h1}
\end{align}
Here $y(x)$ is the solution of the Cauchy problem~\eqref{eq:*1}, and $y_1(x)$ and~$y_2(x)$  are
the solutions of the corresponding auxiliary Cauchy problems:
\begin{align}
y^{\prime}_x&=f(0,y)\quad \ \ (x>0),\quad \ y(0)=y_0;\label{*12.1.10.32h2}\\
y^{\prime}_x&=f(I_1,y)\quad \ (x>0),\quad \ y(0)=y_0.\label{*12.1.10.32h2a}
\end{align}
The solutions $y_1(x)$ and~$y_2(x)$ can be represented implicitly as follows:
\begin{align}
x=\int^y_{y_0}\frac{d\xi}{f(0,\xi)},\quad \
x=\int^y_{y_0}\frac{d\xi}{f(I_1,\xi)}.
\label{*12.1.10.32h2*}
\end{align}
For the critical value~$x_*$, the two-sided estimate
\begin{align}
I_2\le x_*\le I_1
\label{*12.1.10.32h3}
\end{align}
is valid.

\textsl{Case $2^\circ$.} Let $f_x\le 0$.
Suppose that the integrals~$I_1$ and~$I_2$ in~\eqref{*12.1.10.32hz} exist and are finite.
Suppose also that the conditions,
\begin{align}
f(x,y)>0,\quad f_x(x,y)\le 0 \quad \text{for all}\quad 0\le x\le I_2, \ \ y\ge y_0>0,
\label{*12.1.10.32g**}
\end{align}
are satisfied. Then the following inequalities are valid:
\begin{align}
f(I_1,y)\le f(x,y)\le  f(0,y) \quad \text{for}\quad 0\le x\le I_2
\label{*12.1.10.32h**}
\end{align}
and
\begin{align}
y_2(x)\le y(x)\le y_1(x) \quad \text{for}\quad 0\le x\le I_1\le I_2,
\label{*12.1.10.32h1**}
\end{align}
where $y(x)$ is the solution of the Cauchy problem~\eqref{eq:*1}, and $y_1(x)$ and~$y_2(x)$
are the solutions of the corresponding auxiliary Cauchy problems \eqref{*12.1.10.32h2} and \eqref{*12.1.10.32h2a}.
The last two solutions can be represented implicitly~\eqref{*12.1.10.32h2*}.
For the critical value~$x_*$, the two-sided estimate
\begin{align}
I_1\le x_*\le I_2
\label{*12.1.10.32h3**}
\end{align}
is valid.

\begin{example}
We consider the Cauchy problem for the Riccati equation
\begin{align}
y^{\prime}_x=y^2+h(x)\quad (x>0);\quad \ y(0)=a>0.
\label{*12.1.10.32f.R*}
\end{align}

Let us consider the two cases.

\textsl{Case $1^\circ$.}  Let $h(x)\ge 0$ and $h'_x(x)\ge 0$.
In this case, the first auxiliary
Cauchy problem~\eqref{*12.1.10.32h2} is written as follows:
\begin{align}
y^{\prime}_x=y^2+h(0)\quad \ (x>0),\quad \ y(0)=a.
\label{*12.1.10.32h4}
\end{align}
The exact solution of the problem~\eqref{*12.1.10.32h4} admits an implicit form
of representation with the help of the first relation~\eqref{*12.1.10.32h2*}
for $y_0=a$ and $f(0,y)=y^2+h(0)$. After elementary calculations and transformations,
this solution can be written in the explicit form
\begin{align}
y=\sqrt b\,\frac{a\cos(\sqrt b\,x)+\sqrt b \sin(\sqrt b\,x)}{\sqrt b \cos(\sqrt b\,x)-a\sin(\sqrt b\,x)},\quad \ b=h(0).
\label{*12.1.10.32h5}
\end{align}
The singular point of this solution, $I_1$, which is the zero of the denominator and is equal to the improper first integral
in~\eqref{*12.1.10.32hz} for $y_0=a$ and $f(0,y)=y^2+h(0)$, is defined by the formula
$$
I_1=\frac 1{\sqrt b}\arctan\frac{\sqrt b}a,\quad \ b=h(0).
$$

The solution of the second auxiliary Cauchy problem~\eqref{*12.1.10.32h2a} is given by the
formula \eqref{*12.1.10.32h5}, in which  $h(0)$ must be replaced by~$h(I_1)$.
As a result, we obtain the two-sided estimate of the critical value~$x_*$:
\begin{equation}
\begin{gathered}
I_2\le x_*\le I_1,\\
I_1=\frac 1{\sqrt{h(0)}}\arctan\frac{\sqrt {h(0)}}a,\quad
I_2=\frac 1{\sqrt{h(I_1)}}\arctan\frac{\sqrt {h(I_1)}}a.
\end{gathered}
\label{*12.1.10.32h7}
\end{equation}

\begin{figure}
\centering
{\includegraphics[scale=0.34]{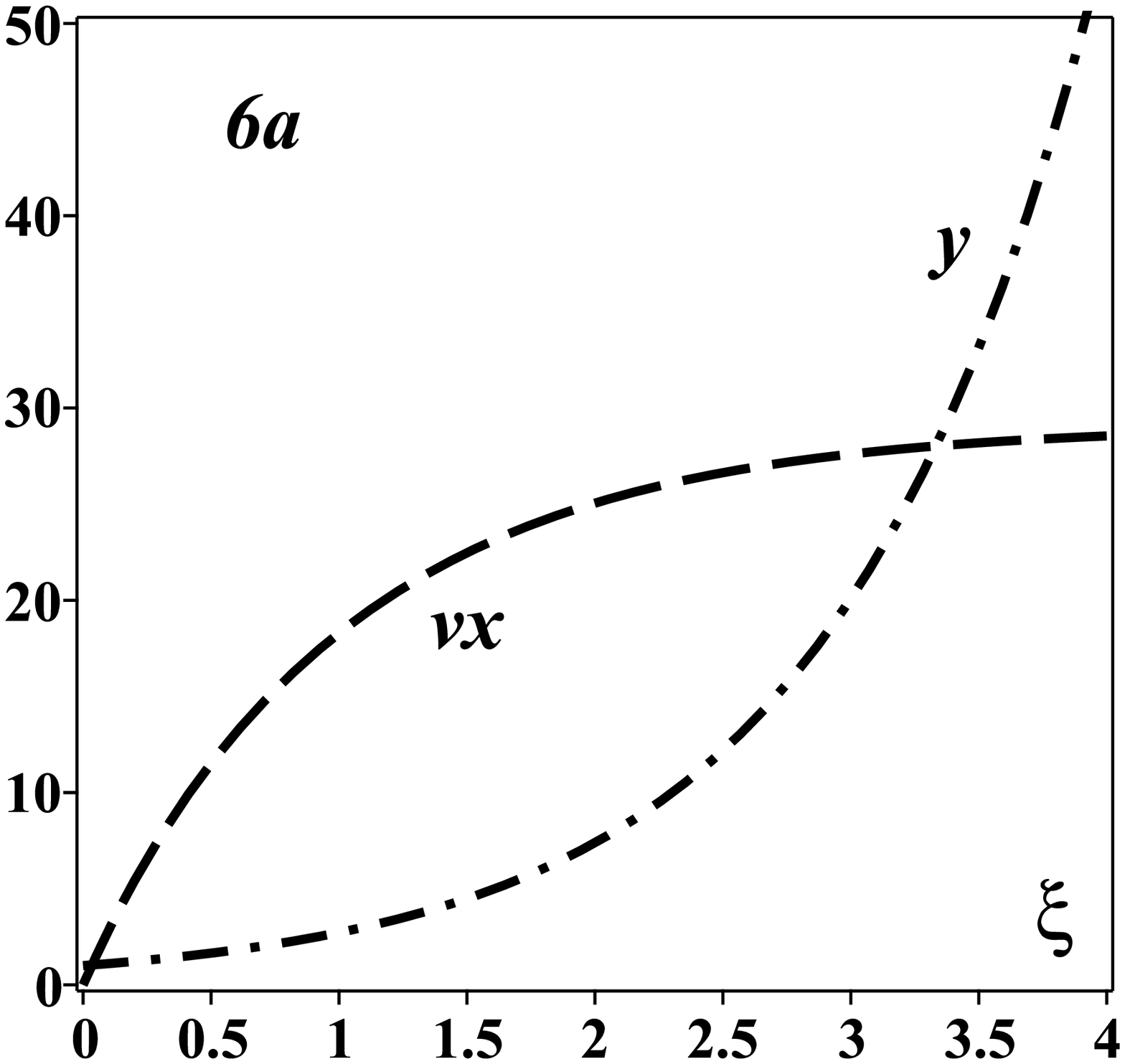}\ \includegraphics[scale=0.34]{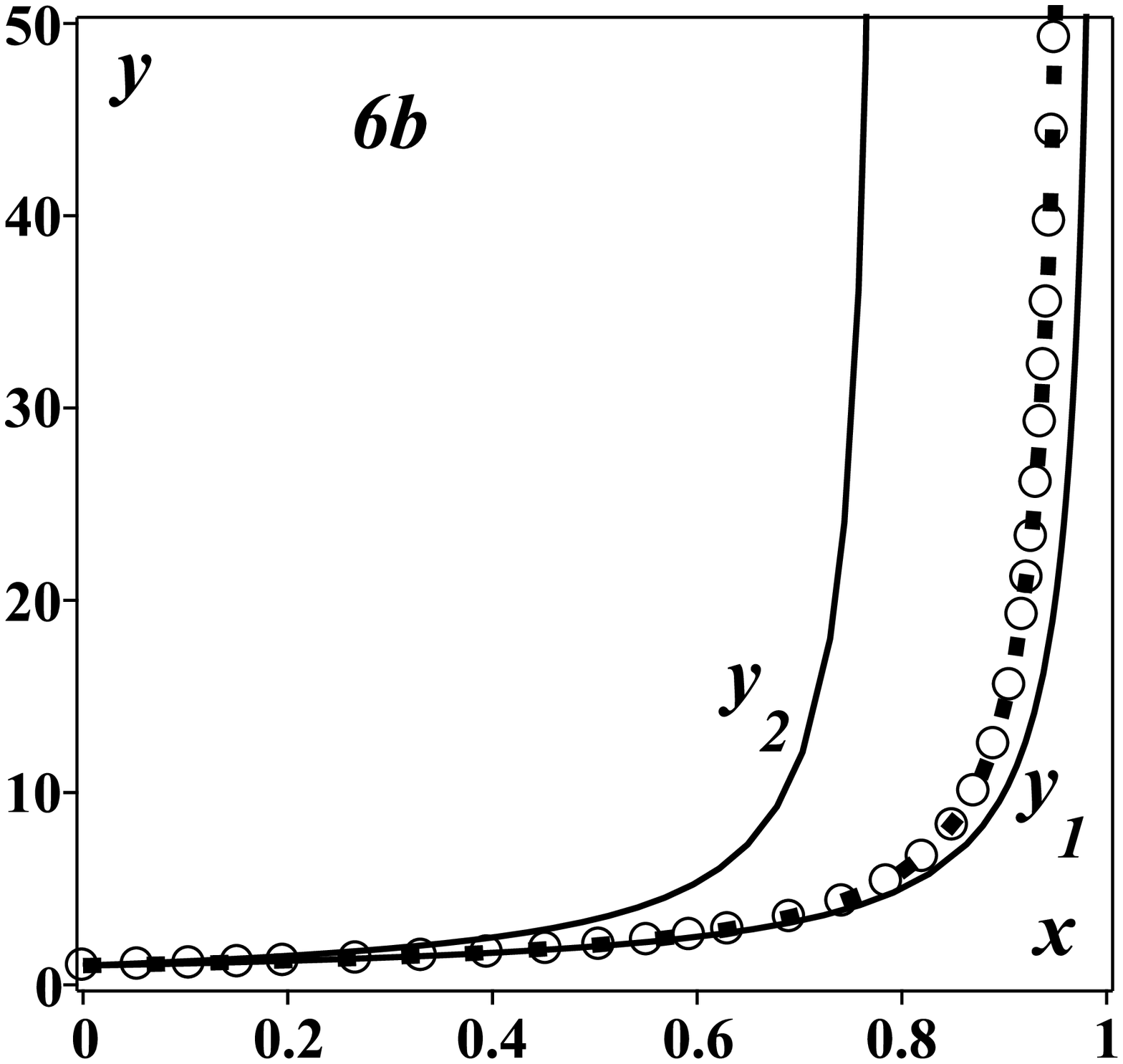}}
\caption{{\it 6a}\dash the dependences $x=x(\xi)$ and $y=y(\xi)$ obtained numerically after the transformation
of the Cauchy problem for one equation~\eqref{*12.1.10.32f.R*} for $a=1$, $h(x)=x^2$ to the problem for
the system of equations~\eqref{eq:02p} for $f=y^2+x^2$, $g=f/y$ ($\nu=30$);
{\it 6b}\dash the exact solution of the problem~\eqref{*12.1.10.32f.R*} (points),
the numerical solution of this problem (circles), and the majorizing functions $y_1(x)$ and $y_2(x)$ (solid lines).}
 \label{fig:Fig4}
\end{figure}

In particular, setting $a=1$, $h(x)=x^m$ and $m>0$ in~\eqref{*12.1.10.32f.R*},
we find that $I_1=1$ and $I_2=\arctan 1$.
Substituting these values into~\eqref{*12.1.10.32h7}, we obtain the two-sided estimate $0.785\le x_*\le 1$
for the critical value~$x_*$.

In Fig.~\ref{fig:Fig4} we present the results of the numerical solution of the Cauchy problem~\eqref{*12.1.10.32f.R*}
for $a=1$ and $h(x)=x^2$  in parametric form, as well as a comparison of the
numerical and exact solutions of this problem (the latter is expressed in terms of the Bessel functions and is omitted here),
and also the majorizing functions $y_1=1/(1-x)$ and $y_2=y_2(x)$,
which are the solutions of the auxiliary Cauchy problems~\eqref{*12.1.10.32h2} and~\eqref{*12.1.10.32h2a}
(the solution of the Cauchy problem under consideration is located between these functions).
The function $y_2(x)$ is determined by the formula~\eqref{*12.1.10.32h5}, in which the parameter~$b$ must be replaced by~$1/a$.

\medskip
We note that if $h(x)=\text{const}>0$, then the inequalities~\eqref{*12.1.10.32h7} give the exact result
\text{$x_*=I_1=I_2$}.

\textsl{Case $2^\circ$.}  Let $h(x)\ge 0$ and $h'_x(x)\le 0$.
In this case, the solution of the first auxiliary Cauchy problem~\eqref{*12.1.10.32h4}
is also given by the formula~\eqref{*12.1.10.32h5}, and the solution of the second
auxiliary Cauchy problem is obtained from~\eqref{*12.1.10.32h5} by a formal replacement of $h(0)$ by~$h(I_1)$.
As a result, we obtain the two-sided estimate for the critical value~$x_*$:
$$
I_1\le x_*\le I_2,
$$
where the integrals $I_1$ and $I_2$ are determined by the formulas~\eqref{*12.1.10.32h7}.
\end{example}

\begin{remark}
It should be noted that in Case $2^\circ$ it does not matter how
the function $f(x,y)$ and its derivative $f_x(x,y)$ behave for $x>I_2$; in particular,
the right-hand side of the equation~\eqref{eq:*1} can be negative for $x>I_2$.
\end{remark}

\begin{example}
To illustrate what was said in Remark~9, we consider the test Cauchy problem
\begin{align}
y^{\prime}_x=(2-x)y^2\quad (x>0);\quad \ y(0)=1,
\label{eq:500}
\end{align}
which corresponds to \textsl{Case $2^\circ$}, where $f(x,y)<0$ for $x>2$.

Calculating the integrals~\eqref{*12.1.10.32hz}, we have $I_1=\frac12$ and $I_2=\frac23$.
Substituting these values into~\eqref{*12.1.10.32h3**}, we obtain the
two-sided estimate for the singular point
$$
\tfrac12<x_*<\tfrac23.
$$

The exact solution of the problem~\eqref{eq:500} is given by the formula
\begin{align}
y=\frac 2{x^2-4x+2}.
\label{eq:501}
\end{align}
The zero of the denominator, equal to $x_*=2-\sqrt 2\approx 0.5858$,
determines the singular point of the solution (first-order pole).
\end{example}

\section{Problems for second-order equations. Differential transformations}\label{ss:6}

\subsection{Solution method based on introducing a differential variable}\label{ss:6.1}

The Cauchy problem for the second-order differential equation has the form
\begin{align}
&y^{\prime\prime}_{xx}=f(x,y,y^{\prime}_x)\quad (x>x_0);\label{eq:03}\\
&y(x_0)=y_0,\quad \ y^{\prime}_x(x_0)=y_1.\label{eq:03a}
\end{align}

We note that the exact solutions of equations of the form~\eqref{eq:03}, which can be used
for the formulation of test problems with blow-up solutions, can be found in
\cite{pol2003,mur,kam,kudr}.

Let $f(x,y,u)>0$ if $y>y_0\ge 0$ and $u>y_1\ge 0$, and the function~$f$
increases quite rapidly as $y\to\infty$ (for example, if $f$ does not depend on~$y^{\prime}_x$, then
$\ds\lim_{y\to\infty}f/y=\infty$).

First, as in Section~2.1, we represent the ODE~\eqref{eq:03}
as an equivalent system of differential-algebraic equations
\begin{align}
y^{\prime}_x=t,\quad \ y^{\prime\prime}_{xx}=f(x,y,t),
\label{eq:03b}
\end{align}
where $y=y(x)$ and $t=t(x)$ are the unknown functions.

Taking into account~\eqref{eq:03b}, we derive a standard system of ODEs
for the functions $y=y(t)$ and $x=x(t)$. To do this, differentiating
the first equation of the system~\eqref{eq:03b}  with respect to~$t$, we obtain  $(y^{\prime}_x)'_t=1$.

Taking into account the relations $y^{\prime}_t=tx^{\prime}_t$ (it follows from the first equation~\eqref{eq:03b}) and
$(y^{\prime}_x)'_t=y^{\prime\prime}_{xx}/t'_x=x^{\prime}_ty^{\prime\prime}_{xx}$, we have
\begin{equation}
x^{\prime}_ty^{\prime\prime}_{xx}=1.
\label{eq:03c}
\end{equation}
Eliminating here the second derivative $y^{\prime\prime}_{xx}$ by using the second equation~\eqref{eq:03b},
we arrive at the first-order equation
\begin{equation}
x^{\prime}_t=\frac 1{f(x,y,t)}.
\label{eq:03d}
\end{equation}
Considering further the relation $y^{\prime}_t=tx^{\prime}_t$, we transform~\eqref{eq:03d} to the form
\begin{equation}
y^{\prime}_t=\frac t{f(x,y,t)}.
\label{eq:03e}
\end{equation}
Equations \eqref{eq:03d} and \eqref{eq:03e} represent a system of coupled first-order differential equations
for the unknown functions $x=x(t)$ and $y=y(t)$.
The system \eqref{eq:03d}--\eqref{eq:03e} should be supplemented by the initial conditions
\begin{align}
x(t_0)=x_0,\quad y(t_0)=y_0,\quad t_0=y_1,
\label{eq:03f}
\end{align}
which are derived from~\eqref{eq:03a} and the first equation~\eqref{eq:03b}.

The Cauchy problem \eqref{eq:03d}--\eqref{eq:03f} has a solution without
blow-up singularities and can be integrated by applying the standard fixed-step numerical methods (see, for example,
\cite{but,fox87,dormand96,sch,shampine94,asch,korn,shing2009,grif2010}).
\medskip

\begin{remark}
Systems of differential-algebraic equations \eqref{eq:02b} and \eqref{eq:03b}
are particular cases of parametrically defined nonlinear  differential equations, which are considered
in~\cite{pol2016,pol2017}.
In~\cite{pol2017}, the general solutions of several parametrically defined ODEs were
constructed via differential transformations based on introducing a new differential independent
variable  $t=y^{\prime}_x$.
\end{remark}

\subsection{Test problems and numerical solutions}\label{ss:6.2}

\begin{example}
We consider a test Cauchy problem for the second-order nonlinear ODE
\begin{align}
y^{\prime\prime}_{xx}&=b\gamma y^{\gamma-1}y^\prime_x\quad \ \ (x>0);\quad y(0)=a,\quad y^{\prime}_x(0)=a^\gamma b,
\label{eq:04}
\end{align}
which is obtained by differentiating equation~\eqref{eq:xx02}.
For  $a>0$, $b>0$, and $\gamma>1$,
the exact solution of this problem is defined  by the formula~\eqref{eq:xx04}.

Introducing a new variable $t=y^{\prime}_x$ in~\eqref{eq:04},
we obtain the Cauchy problem, which exactly coincides with the problem~\eqref{eq:02h}.
The exact solution of this problem is determined by the formulas~\eqref{eq:02j}.
\end{example}

\begin{figure}
\centering
{\includegraphics[scale=0.34]{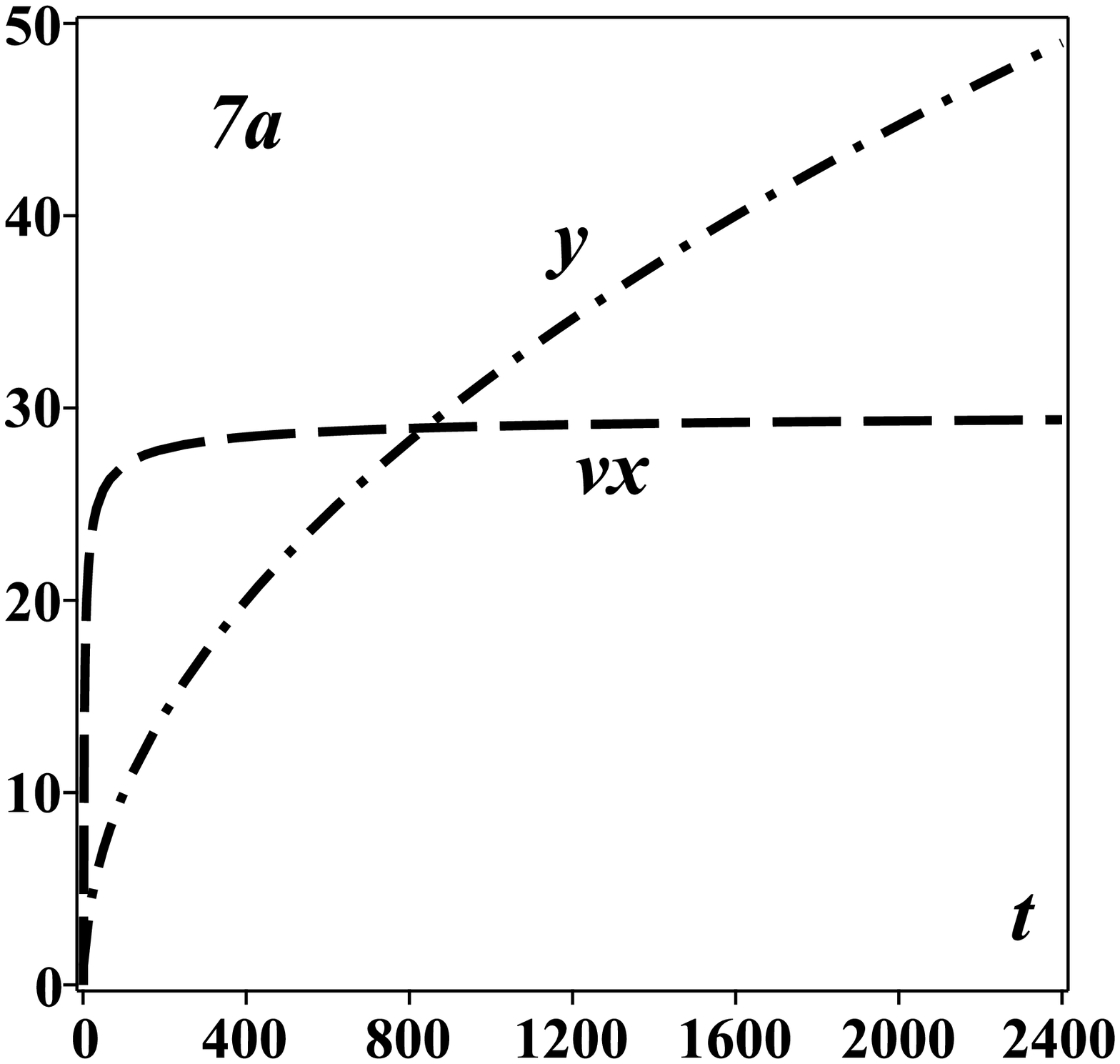}\ \includegraphics[scale=0.34]{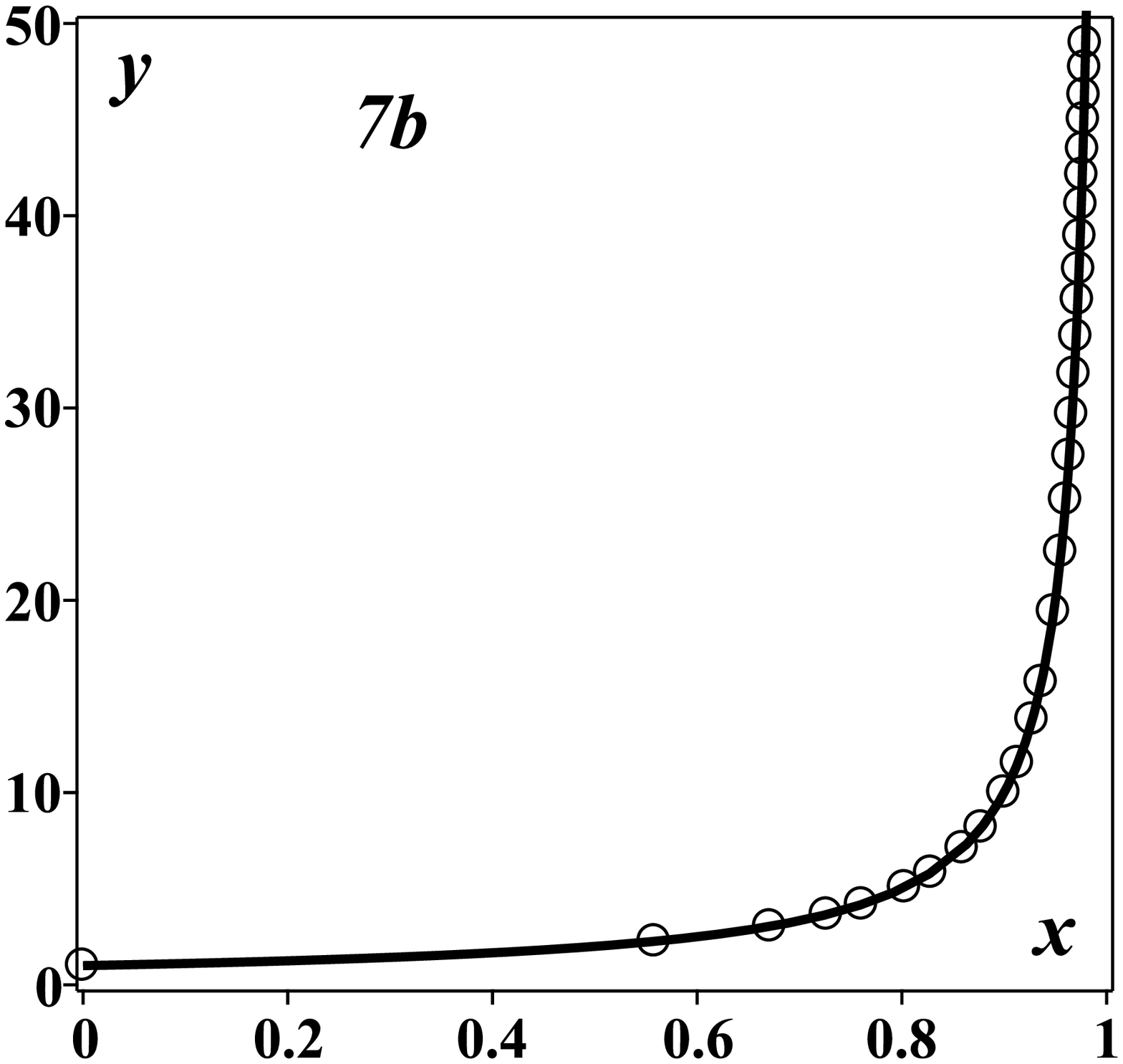}}
\caption{{\it 7a}\dash the dependences $x=x(t)$ and $y=y(t)$ obtained by numerical solution of the problem~\eqref{eq:04b}
for $a=b=1$ and $\gamma=2$ ($\nu=30$); {\it 7b}\dash exact solution~\eqref{eq:02g} (solid line) and numerical solution of problem~\eqref{eq:04b}
(circles).}
\label{fig:Fig5}
\end{figure}

\begin{example}
Let us now consider another Cauchy problem
\begin{align}
y^{\prime\prime}_{xx}&=b^2\gamma y^{2\gamma-1}\quad \ \ (x>0);\quad y(0)=a,\quad y^{\prime}_x(0)=a^\gamma b,
\label{eq:04a}
\end{align}
which is obtained by excluding the first derivative from the equations
\eqref{eq:xx02} and~\eqref{eq:04}
(we recall that the second equation is a consequence of the first equation).
The exact solution of the problem~\eqref{eq:04a} is determined by the formula~\eqref{eq:xx04}.

Introducing a new variable $t=y^{\prime}_x$, we transform~\eqref{eq:04a} to the Cauchy problem
for the system of the first-order ODEs
\begin{equation}
\begin{gathered}
x^{\prime}_t=\frac 1{b^2\gamma y^{2\gamma-1}},\quad \ y^{\prime}_t=\frac t{b^2\gamma y^{2\gamma-1}}\quad \ (t>t_0);\\
x(t_0)=0,\quad y(t_0)=a,\quad t_0=a^\gamma b,
\end{gathered}
\label{eq:04b}
\end{equation}
which is a particular case of the problem \eqref{eq:03d}--\eqref{eq:03f}
with $f=b^2\gamma y^{2\gamma-1}$, $x_0=0$, and $y_0=a$.
The exact solution of the problem~\eqref{eq:04b} is given by formulas~\eqref{eq:02j}.
\medskip

Figure~\ref{fig:Fig5} shows a comparison of the exact solution~\eqref{eq:02g} of the Cauchy problem
for one equation~\eqref{eq:02f} with the numerical solution of the related problem for the system of equations~\eqref{eq:04b} for $a=b=1$ and $\gamma=2$,
obtained by applying the Runge--Kutta method of the fourth-order of approximation.

The function $x(t)$ slowly tends to the asymptotic value~$x_*$. Therefore
to accelerate this process in the system~\eqref{eq:04b} is useful additionally to  make
the exponential-type substitution~\eqref{eq:02k}.
\end{example}

For completeness of the picture, we also give an example of a blow-up problem
whose solution has a logarithmic singularity.

\begin{example}
An exact solution of the Cauchy problem with exponential nonlinearity
\begin{align}
y^{\prime\prime}_{xx}=e^{2y}\quad (x>0);\qquad y(0)=0,\quad y^{\prime}_x(0)=1,
\label{12.1.10.31***}
\end{align}
has the form
\begin{align}
y=\ln\Bl(\frac{1}{1-x}\Br)=-\ln(1-x).
\label{12.1.10.32***}
\end{align}
This solution has a logarithmic singularity at the point  $x_*=1$ and
does not exist for $x>x_*$.

Introducing the differential variable $t=y^{\prime}_x$, we transform
the problem~\eqref{12.1.10.31***} to the following Cauchy problem for a system of equations:
\begin{equation}
\begin{gathered}
x^{\prime}_t=e^{-2y},\quad \ y^{\prime}_t=te^{-2y}\quad \ (t>1);\\
x(1)=0,\quad \ y(1)=0\quad \ (t_0=1),
\end{gathered}
\label{eq:505}
\end{equation}
which is a particular case of the system \eqref{eq:03d}--\eqref{eq:03e}.
The exact analytical solution of the problem~\eqref{eq:505} is determined by the formulas
\begin{align}
x=1-\frac 1t,\quad \ y=\ln t\qquad (t\ge 1),
\notag
\end{align}
which do not have singularities; the function $x=x(t)$ increases monotonically
with  $t>1$ and tends to its limiting value $\ds x_*=\lim_{t\to\infty}x(t)=1$, and
the function $y=y(t)$ is unlimited and increases monotonically with respect to the logarithmic law.

The function $x(t)$ slowly tends to the asymptotic value~$x_*$. Therefore
to accelerate this process in the system~\eqref{eq:505} is useful additionally to  make
the exponential-type substitution~\eqref{eq:02k}.
\end{example}

\section{Problems for second-order equations. Nonlocal transformations and\\ differential constraints}\label{ss:7}

\subsection{Solution method based on introducing a non-local variable}\label{ss:7.1}

First, we represent the equation~\eqref{eq:03} as the equivalent system of two equations
\begin{align}
y^{\prime}_x=t,\quad \ t'_x=f(x,y,t),
\notag
\end{align}
and then we introduce a non-local variable of the form \cite{pol_shi2017}:
\begin{align}
\xi=\int^x_{x_0}g(x,y,t)\,dx,\quad \ y=y(x),\quad t=t(x).
\label{eq:05}
\end{align}

As a result, the Cauchy problem \eqref{eq:03}--\eqref{eq:03a} is transformed to the following problem for
the autonomous system of three equations:
\begin{equation}
\begin{gathered}
x^{\prime}_\xi=\frac 1{g(x,y,t)},\quad \ y^{\prime}_\xi=\frac t{g(x,y,t)},\quad \ t'_\xi=\frac {f(x,y,t)}{g(x,y,t)}\quad \ (\xi>0);\\
x(0)=x_0,\quad \ y(0)=y_0,\quad \ t(0)=y_1.
\end{gathered}
\label{eq:06}
\end{equation}

For a suitable choice of the function $g=g(x,y,t)$
(not very restrictive conditions of the type~\eqref{eq:02q} must be imposed on it),
we obtain the Cauchy problem~\eqref{eq:06}, the solution of which will not have blow-up singularities;
therefore this problem can be integrated by applying the standard fixed-step numerical methods \cite{but,fox87,dormand96,sch,shampine94,asch,korn,shing2009,grif2010}.

\medskip
Let us consider various possibilities for choosing the function~$g$ in the system~\eqref{eq:06}.

\begin{description}

\item{$1^\circ$.} The special case $g=t$ is equivalent to the hodograph transformation with
an additional shift of the dependent variable, which gives $\xi=y-y_0$.

\item{$2^\circ$.} We can take $g=\bl(c+|t|^s+|f|^s\br)^{1/s}$ for $c\ge 0$ and $s>0$. The case $c=1$ and $s=2$
corresponds to the method of the arc-length transformation~\cite{mor1979}.

\item{$3^\circ$.} Setting $g=f$, and then integrating the third equation~\eqref{eq:06},
we obtain the problem \eqref{eq:03d}--\eqref{eq:03f} in which the variable  $t=\xi+y_1$.
Therefore the method based on the non-local transformation~\eqref{eq:05} generalizes the method
based on the differential transformation, which is described in Section~6.1.


\item{$4^\circ$.} We can take $g=t/y$ (or $g=kt/y$, where $k>0$ is a numerical parameter that can be varied).
In this case, the system~\eqref{eq:06} is much simplified, since the second equation is directly integrated,
and taking into account the second initial condition, we obtain $y=y_0e^\xi$.
As a result, there remains a system of two equations for the determination of the functions $x=x(\xi)$ and $t=t(\xi)$.
Taking into account the relation $t=y^\prime_x$, we also have
$$
\xi=\int^x_{x_0}\frac{t}y\,dx=\int^x_{x_0}\frac{y^{\prime}_{x}}{y}\,dx=\ln\frac{y}{y_0}.
$$
Therefore, the non-local transformation~\eqref{eq:05} with $g=t/y$ and the subsequent transition to the system~\eqref{eq:06}
is equivalent to a point transformation $\xi=\ln(y/y_0)$, $z=x$, which is a combination of
two more simple point transformations: 1) the transformation $\bar x=x$, $\bar y=\ln(y/y_0)$
and 2) the hodograph transformation $\xi=\bar y$, $z=\bar x$, where $z=z(\xi)$.

\item{$5^\circ$.} Also, we can take $g=f/t$ (or $g=kf/t$, where $k>0$ is a free numerical parameter).
In this case, the system~\eqref{eq:06} is also much simplified, since the third equation is directly integrated,
and taking into account the third initial condition, we obtain $t=y_1e^\xi$.
As a result, there remains a system of two equations for the determination of the functions $x=x(\xi)$ and $y=y(\xi)$.
For the non-local transformation~\eqref{eq:05} with $g=f/t$, the new independent variable is expressed
in terms of the derivative by the formula $\xi=\ln(y^{\prime}_x/y_1)$
(that is, this transformation coincides with the modified differential transformation, see Section~2.3).


\end{description}

\begin{remark}
The transformations corresponding to the last two cases, $4^\circ$ and $5^\circ$,
will be called special exp-type transformations, they lead
to the solutions, in which the variable~$x$ tends exponentially rapidly to a blow-up point~$x_*$.
\end{remark}

\begin{remark}
From Items~$1^\circ$, $2^\circ$, and $3^\circ$
it follows that the method based on the hodograph transformation, the method of the arc-length transformation,
and the method based on the differential transformation are particular cases of the non-local transformation
of the general form~\eqref{eq:05}, which leads to the Cauchy problem for the system of equations~\eqref{eq:06}.
\end{remark}

\begin{remark}
It is not necessary to calculate the integrals~\eqref{eq:05} (or~\eqref{*})
when using non-local transformations.
\end{remark}

\subsection{Test problems and numerical solutions}\label{ss:7.2}

\begin{example}
For the test problem~\eqref{eq:04a}, in which $f=b^2\gamma y^{2\gamma-1}$, we set
$g=t/y$ (see Item~$4^\circ$ in Section~7.1). Substituting these functions
into~\eqref{eq:06}, we arrive at the Cauchy problem
\begin{equation}
\begin{gathered}
x^{\prime}_\xi=\frac yt,\quad \ y^{\prime}_\xi=y,\quad \ t'_\xi=\frac{b^2\gamma y^{2\gamma}}t\quad \ (\xi>0);\\
x(0)=0,\quad \ y(0)=a,\quad \ t(0)=a^\gamma b.
\end{gathered}
\label{eq:xx}
\end{equation}
The exact solution of this problem in parametric form is determined by the formulas
\begin{equation}
x=\frac 1{a^{\gamma-1}b(\gamma-1)}\bl[1-e^{-(\gamma-1)\xi}\br],\quad \ y=ae^\xi,\quad \ t=a^\gamma be^{\gamma\xi}.
\label{eq:xy}
\end{equation}
\noindent
It can be seen that the required function $x=x(\xi)$ exponentially tends to
the asymptotic value $\ds x_*=\frac 1{a^{\gamma-1}b(\gamma-1)}$ as $\xi\to \infty$.
\medskip

\begin{figure}
\centering
{\includegraphics[scale=0.34]{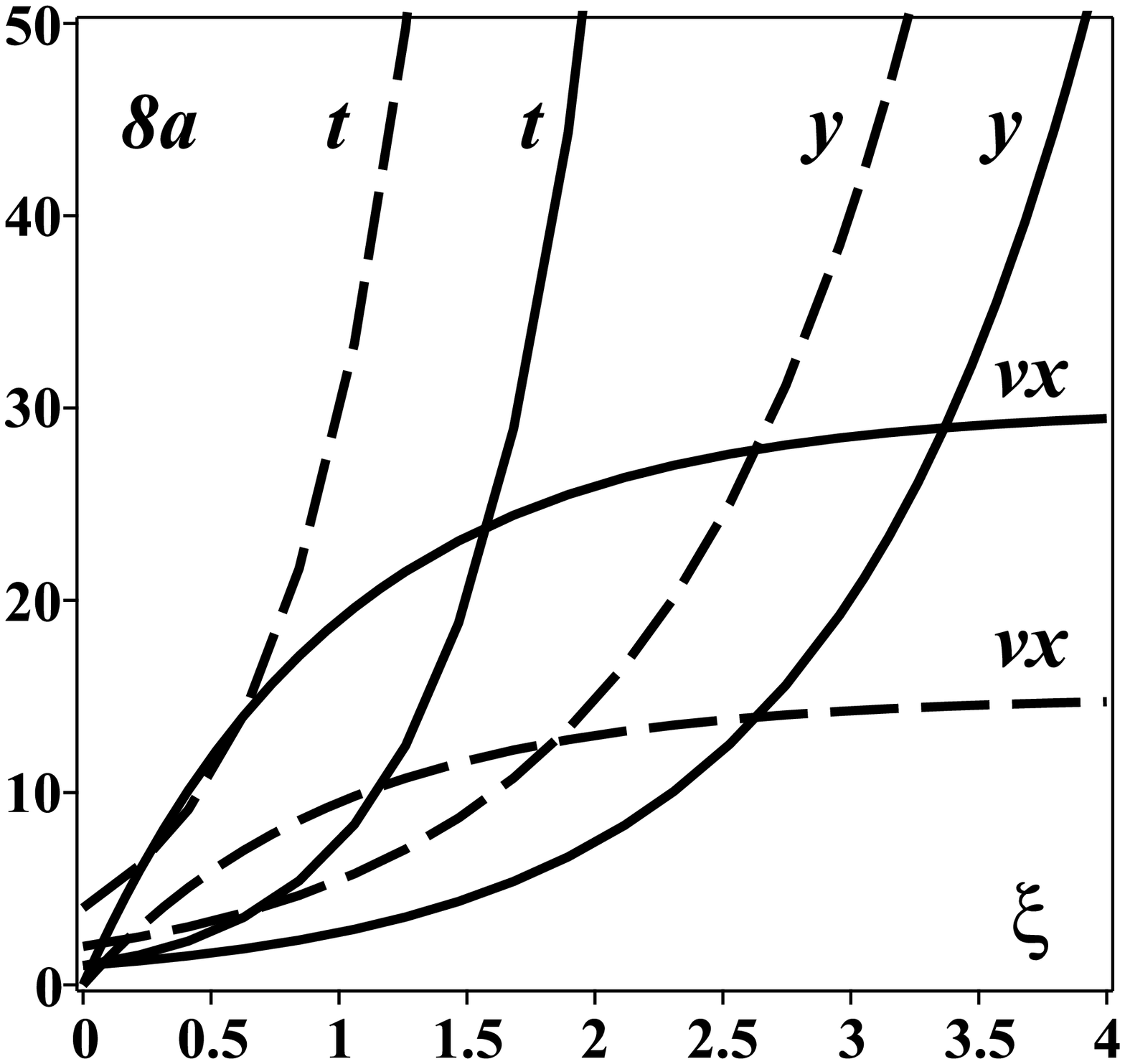}\ \includegraphics[scale=0.34]{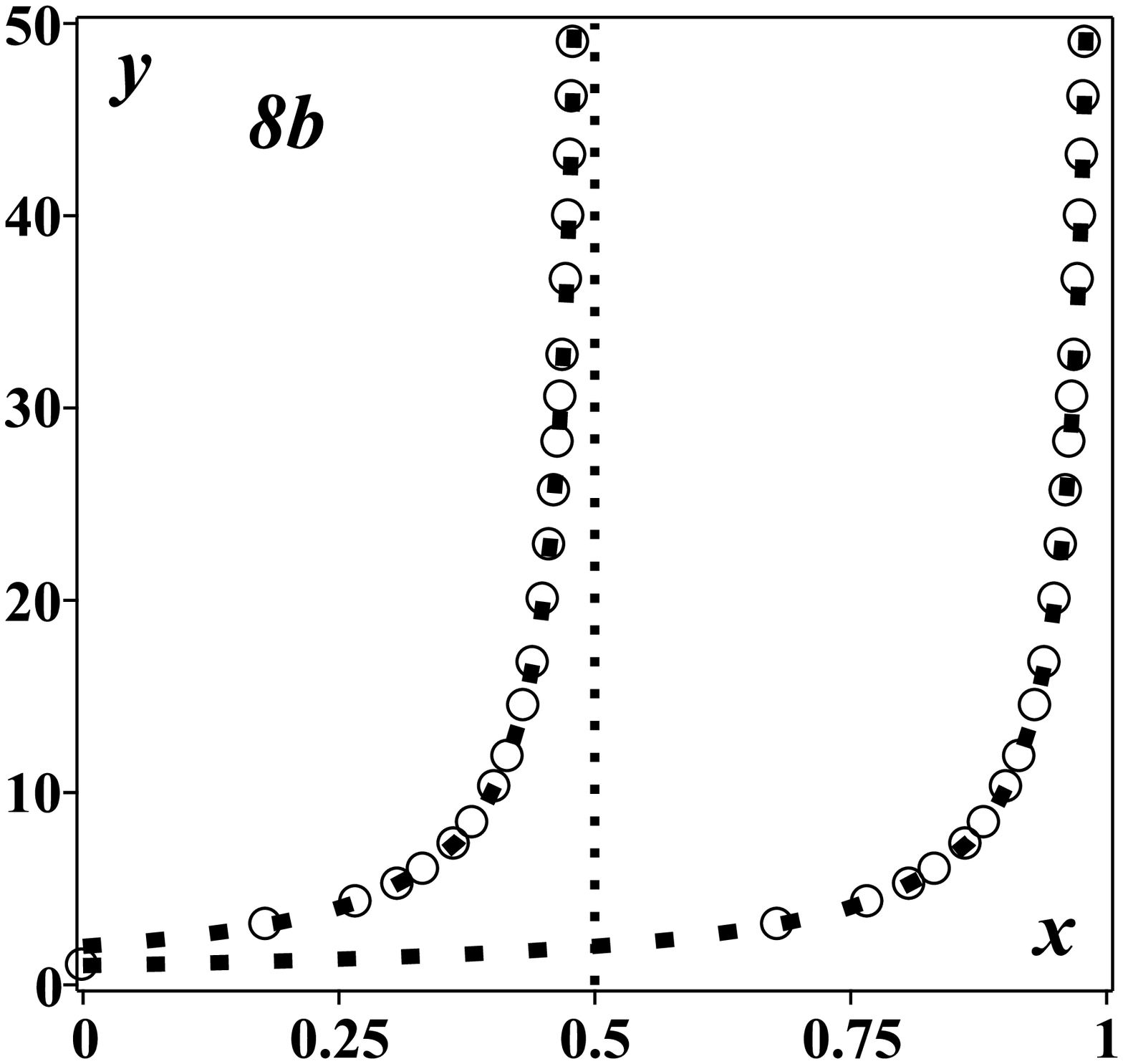}}
\caption{{\it 8a}\dash the dependences $x=x(\xi)$, $y=y(\xi)$, $t=t(\xi)$, obtained by numerical solution
of the problem~\eqref{eq:xx} for $b=1$, $\gamma=2$ with $a=1$ (solid lines) and $a=2$ (dashed lines) ($\nu=30$);
{\it 8b}\dash numerical solutions of the problems~\eqref{eq:04b} for $b=1$, $\gamma=2$ (circles)
and~\eqref{eq:xx} (points); for left curves $a=2$ and for right curves $a=1$.}
\label{fig:Fig6}
\end{figure}

The numerical solutions of the problems~\eqref{eq:04b} and~\eqref{eq:xx} for $b=1$, $\gamma=2$,
obtained by the Runge--Kutta method of the fourth-order of approximation,
are shown in Fig.~\ref{fig:Fig6} for $a=1$ and $a=2$. For a fixed step of integration,
equal to~$0.2$, the maximum difference between the exact solution~\eqref{eq:02g}
and the numerical solution of the related problem~\eqref{eq:xx} is $0.0045\%$.
For larger stepsizes, equal to $0.4$ and~$0.6$, the maximum error in the numerical solutions is $0.061\%$ and~$0.24\%$, respectively.
It can be seen that the numerical solutions of the problems~\eqref{eq:04b} and~\eqref{eq:xx}
are in a good agreement, but the rates of their approximation to the required asymptote $x=x_*$
are significantly different. For example, for the system~\eqref{eq:04b},
in order to obtain a good approximation to the asymptote, it is required to consider
the interval $t\in[1,2400]$, and for the system~\eqref{eq:xx} it suffices to take
$\xi\in[0,4]$. Therefore, it should expect that the method based on the use of the system~\eqref{eq:06} with $g=t/y$
is much more efficient than the method based on the differential transformation.

For comparison, similar calculations were performed using Maple (2016), and applying
the method based on the hodograph transformation (see Section~7.1, Item~$1^\circ$)
and the method of the arc-length transformation (see Section~7.1, Item~$2^\circ$ for $c=1$ and $s=2$).
In order to obtain a good approximation to the asymptote, applying the method based on the hodograph transformation,
it is required to consider the interval $\xi\in[0,49]$,
while using the method of the arc-length transformation leads to a significantly larger interval
$\xi\in[0,2500]$.
To control a numerical integration process, the calculations were carried out
with the aid of two other most important and powerful
mathematical software packages: Mathematica~(11), and MATLAB (2016a).
It was found that the method based on the use of the system~\eqref{eq:xx}
with $g=t/y$ is essentially
more efficient than the method based on the hodograph transformation
and the method of the arc-length transformation.
\end{example}

\begin{example}
For the test problem~\eqref{eq:04a}, in which $f=b^2\gamma y^{2\gamma-1}$, we set $g=f/t$
(see Item~$5^\circ$ in Section~7.1). Substituting these functions into~\eqref{eq:06},
we arrive at the Cauchy problem
\begin{equation}
\begin{gathered}
x^{\prime}_\xi=\frac t{b^2\gamma y^{2\gamma-1}},\quad \ y^{\prime}_\xi=\frac{t^2}{b^2\gamma y^{2\gamma-1}},\quad \ t'_\xi=t\quad \ (\xi>0);\\
x(0)=0,\quad \ y(0)=a,\quad \ t(0)=a^\gamma b.
\end{gathered}
\label{eq:xx**}
\end{equation}
The exact solution of this problem in parametric form is
\begin{equation}
x=\frac {1}{a^{\gamma-1}b(\gamma-1)}\Bl[1-e^{-(\gamma-1)\xi/\gamma}\Br],\quad \ y=ae^{\xi/\gamma},\quad \ t=a^\gamma be^{\xi}.
\label{eq:xy**}
\end{equation}
\noindent

The required value $x=x(\xi)$ tends exponentially to the asymptotic value $\ds x_*=\frac {1}{a^{\gamma-1}b(\gamma-1)}$ as $\xi\to \infty$.
However, in comparison with the method applied in Example~19,
in this case the rate of approximation of the parametric solution to the asymptote is less
(which is not important for application of the standard numerical methods for solving similar problems).
Note that the solution~\eqref{eq:xy**} coincides with~\eqref{eq:xy} if we redenote $\xi$ by $\gamma\xi$.
\end{example}

\subsection{Generalizations based on the use of differential constraints}\label{ss:7.3}

The method of numerical integration of the Cauchy problems with blow-up solutions, which based
on introducing a non-local variable, can be generalized if the relation~\eqref{eq:05}
is replaced by the first-order differential constraint
\begin{align}
\xi'_x=g(x,y,t,\xi)
\label{eq:xz**}
\end{align}
with the initial condition $\xi(x=x_0)=\xi_0$.

If we set $\xi_0=0$, then the use of the differential constraint~\eqref{eq:xz**}
leads to the problem~\eqref{eq:06}, where the function $g(x,y,t)$ must be replaced by
$g(x,y,t,\xi)$ in the equations.

Using differential constraints increases the possibilities for numerical analysis of blow-up problems.

In particular, if we choose a differential constraint of the form \eqref{eq:xz**} with
\begin{align}
g(x,y,t,\xi)=\frac t{\varphi(\xi)y+\psi(\xi)},
\label{500}
\end{align}
where  $\varphi(\xi)$ and $\psi(\xi)$ are given functions, then the second equation of the
system~\eqref{eq:06} is reduced to the linear equation for $y=y(\xi)$,
the solution of which is well known. As a result, the considered system,
consisting of three equations, is simplified and reduced to two equations.

If we choose a differential constraint of the form \eqref{eq:xz**} with
\begin{align}
g(x,y,t,\xi)=\frac {f(x,y,t)}{\varphi(\xi)t+\psi(\xi)},
\label{500*}
\end{align}
then the third equation of the system~\eqref{eq:06} is reduced to the linear equation for $t=t(\xi)$.
In this case, the system under consideration also is reduced to two equations.


\begin{example}
For the test Cauchy problem~\eqref{eq:04a} with $b=1$ and $\gamma=2$,
we take the differential constraint~\eqref{eq:xz**} with the function~\eqref{500*},
where $f=2y^3$, $\varphi(\xi)=2(1+2\xi)$, and $\psi(\xi)=0$.
As a result, we arrive at the following problem for the ODE system:
\begin{equation}
\begin{gathered}
x^{\prime}_\xi=\frac{t(1+2\xi)}{y^3},\quad \ y^{\prime}_\xi=\frac{t^2(1+2\xi)}{y^3},\quad \ t^{\prime}_\xi=2t(1+2\xi)\quad \ (\xi>0);\\
x(0)=0,\quad \ y(0)=a,\quad \ t(0)=a^2.
\end{gathered}
\label{501}
\end{equation}
The exact solution of the problem in parametric form is
\begin{equation}
x=\frac 1a\bl(1-e^{-\xi-\xi^2}\br),\quad \ y=ae^{\xi+\xi^2},\quad \
t=a^2e^{2(\xi+\xi^2)}.
\label{eq:xy***}
\end{equation}
It can be seen that the required function $x=x(\xi)$  tends much faster
to the asymptotic value $x_*=1/a$ as $\xi\to \infty$ than in Examples~19 and~20.
\end{example}

\subsection{Comparison of efficiency of various transformations for numerical integration of second-order blow-up ODE problems}\label{ss:7.4}

In Table~2, a comparison of the efficiency of the numerical integration methods,
based on various non-local transformations of the form~\eqref{eq:05} and differential
constraints of the form~\eqref{eq:xz**} is presented by using the example of the test blow-up problem
for the second-order ODE~\eqref{eq:04a} with $a=b=1$ and $\gamma=2$.
The comparison is based on the number of grid points needed to make calculations with the same
maximum error (approximately equal to $0.1$ and~$0.005$).

%



\begin{table}[!ht]
\begin{center}
\begin{tabular}{|l|l|r|r|r|l|l|r|r|r|}
\hline
\multicolumn{5}{|c|}{\vphantom{} {\small Error$_{{\rm max}}, \% =0.1$} }     \\ \hline \hline
{\small Transformation or}     &  \quad {\small Function}  &  {\small Max. interval}  & {\small Stepsize}      & {\small Grid points} \\
{\small differential constraint}  &  \qquad \ \ $g$  &  {\small $\xi_{{\rm max}}\qquad$}  & {\small $h\quad$}  & {\small number}  $N$
\\ \hline
{\small Arc-length}        &  {\small $g{=}\sqrt{1{+}t^2{+}f^2}$}            &   2500.0     &   0.4150                & 6024    \\ \hline
{\small Nonlocal, Item~$2^\circ$}  &  {\small $g{=}1{+}|t|{+}|f|$}           &   2544.0     &   0.7550                & 3369    \\ \hline
{\small Hodograph}         &  $g{=}t$                                        &   49.2       &   0.4510                & 109     \\ \hline
{\small Special exp-type, Item~$5^\circ$} &  $g{=}f/t$                       &   7.807      &   0.2110                & 37      \\ \hline
{\small Diff. constraint, p.c. of \eqref{500*}}  &  $g{=}f/[2t(1+2\xi)]$     &   1.55       &   0.0470                & 33      \\ \hline
{\small Special exp-type, Item~$4^\circ$} &  $g{=}t/y$                       &   3.9        &   0.1300                & 30      \\ \hline
{\small Diff. constraint, p.c. of \eqref{500}}  &  $g{=}t/[2(\xi+1)e^{2\xi+\xi^2}]$ & 1.218 &   0.0435                & 28      \\ \hline
\hline
\multicolumn{5}{|c|}{\vphantom{} {\small Error$_{{\rm max}}, \% =0.005$} }     \\ \hline \hline
{\small Transformation or}     &  \quad Function   &  {\small Max. interval}            & {\small Stepsize}    & {\small Grid points}  \\
{\small differential constraint}  &  \qquad \ \ $g$  &  {\small $\xi_{{\rm max}}\qquad$}  & {\small $h\quad$}  & {\small number}  $N$
\\ \hline
{\small Arc-length}        &  {\small $g{=}\sqrt{1{+}t^2{+}f^2}$}            &   2500.0   &   0.200                & 12500  \\ \hline
{\small Nonlocal, Item~$2^\circ$}  &  {\small $g{=}1{+}|t|{+}|f|$}           &   2544.0   &   0.350                & 7268   \\ \hline
{\small Hodograph}         &  $g{=}t$                                        &   49.0     &   0.125                & 392    \\ \hline
{\small Special exp-type, Item~$5^\circ$} &  $g{=}f/t$                       &   7.821    &   0.099                & 79     \\ \hline
{\small Diff. constraint, p.c. of \eqref{500*}}  &  $g{=}f/[2t(1+2\xi)]$     &   1.55     &   0.021                & 74     \\ \hline
{\small Special exp-type, Item~$4^\circ$} &  $g{=}t/y$                       &   3.9      &   0.060                & 65     \\ \hline
{\small Diff. constraint, p.c. of \eqref{500}} & $g{=}t/[2(\xi+1)e^{2\xi+\xi^2}]$ & 1.220 &   0.020                & 61     \\ \hline
\end{tabular}
\caption{Various types of analytical transformations applied for numerical
integration of the problem~\eqref{eq:04a} for
$a=b=1$ and $\gamma=2$ with a given accuracy (percent errors are $0.1$ and $0.005$
for $\Lambda_\text{m}\le 50$) and their basic parameters  (maximum interval,
stepsize, grid points number). The abbreviation ``p.c.'' stands for ``particular case'' and the notation $f=2y^3$ is used.}
\label{tab:Transfs-ODE2}
\end{center}
\end{table}






It can be seen that the arc-length transformation is the least effective,
since the use of this transformation is associated with a large number of grid points.
In particular, when using the last four transformations, you need 150--200 times less
of a number of grid points.
The hodograph transformation has an intermediate (moderate) efficiency.
The use of the exp-type transformation with $g=t/y$ gives rather good results.

\section{Second-order autonomous equations. Solution of the Cauchy problem. Simple estimates}\label{ss:8}

We consider the Cauchy problem for the second-order autonomous equation
of the general form
\begin{align}
y^{\prime\prime}_{xx}=f(y)\quad (x>0),\quad \ y(0)=a,\quad \ y^{\prime}_x(0)=b.
\label{700}
\end{align}
We assume that $a>0$, $b\ge 0$ and $f(y)>0$ is a continuous function that is defined for all $y\ge a$.

It is not difficult to show that,
the equation~\eqref{700} admits a first integral. As a
result, with allowance for the initial conditions, we arrive at the Cauchy problem for the first-order autonomous equation
\begin{equation}
\begin{gathered}
y^{\prime}_{x}=F(y)\quad (x>0),\quad \ y(0)=a;\\
F(y)=\BL[2\int^y_a f(z)\,dz+b^2\BR]^{\!1/2},
\end{gathered}
\label{701}
\end{equation}
which coincides with the problem~\eqref{12.1.1.1a}, up to obvious modifications in notations.
Therefore, we can use the results of Section~5.1.

The exact solution of the Cauchy problem~\eqref{701} is determined by the formula~\eqref{*12.1.10.32a},
in which the function $f(y)$ should be replaced by~$F(y)$.
In blow-up problems, the critical value~$x_*$ is found by the formula~\eqref{*12.1.10.32b},
where the function $f(y)$ also must be replaced by~$F(y)$.

\medskip
\noindent \textsl{Sufficient criterion of the existence of a blow-up solution.}
Suppose that for some $\kappa>0$ we have the limiting relation
\begin{align}
\lim_{y\to\infty}\frac{F(y)}{y^{1+\kappa}}= s,\quad \ 0<s\le\infty.
\label{702}
\end{align}
Then the solution of the Cauchy problem~\eqref{700}, when the above conditions are satisfied,
is a blow-up solution.

The condition~\eqref{702} is inconvenient, since it contains the function $F(y)$,
which is rather complexly connected with the right-hand side $f(y)$
of the original equation~\eqref{700}.
This condition can be simplified and transformed to a more convenient form:
\begin{align}
\lim_{y\to\infty}\frac{f(y)}{y^{1+\kappa_1}}= s_1,\quad \ 0<s_1\le\infty,
\notag
\end{align}
where $\kappa_1$ is a positive number.

By applying the sufficient criterion, we obtain the following useful result.
\medskip

The Cauchy problem~\eqref{700} for an autonomous equation with power
nonlinearity, $f(y)=cy^\sigma$ ($c>0$), has a blow-up solution if $\sigma>1$.

\section{Blow-up problems for systems of ODEs}\label{ss:9}

\subsection{Method based on non-local transformations}\label{ss:9.1}

We consider the Cauchy problem for a system consisting of~$n$ first-order coupled ODEs
of the general form
\begin{equation}
\frac{dy_m}{dx}=f_m(x,y_1,\dots,y_n),\quad \ \ m=1,\dots, n\quad \ (x>x_0),
\label{eq:08}
\end{equation}
with the initial conditions
\begin{align}
y_m(x_0)=y_{m0},\quad \ \ m=1,\,2,\,\dots,\,n.
\label{eq:08a}
\end{align}

In blow-up problems, the right-hand side of at least one of the equations~\eqref{eq:08}
(after substituting the solution into it) tends to infinity
as $x\to x_*$, where the value~$x_*$ is unknown in advance.

In the general case, the functions $f_m$ may have different signs.
Further, we assume that  $\ds \sum^n_{m=1}|f_m|>0$.

We associate the system~\eqref{eq:08} with the equivalent system of equations consisting of $(n+1)$ equations
\begin{align}
\hskip-1pc \frac{dx}{d\xi}=\frac{1}{g(x,y_1,\dots,y_n)}, \ \ \frac{dy_m}{d\xi}=\frac{f_m(x,y_1,\dots,y_n)}{g(x,y_1,\dots,y_n)}, \ \ m{=}1,\dots, n\ \ (\xi{>}0)
\label{eq:08b}
\end{align}
with the initial conditions
\begin{align}
x(0)=x_0,\quad \ y_m(0)=y_{m0},\quad \ \ m=1,\,2,\,\dots,\,n.
\label{eq:08c}
\end{align}
Here $\xi$ is a non-local variable defined by the formula
\begin{align}
\xi=\int^x_{x_0}g(x,y_1,\dots,y_n)\,dx,\quad \ y_m=y_m(x),\quad \ \ m=1,\dots, n \ \ \ (\xi\ge 0).
\label{eq:08d}
\end{align}

In \eqref{eq:08b}, it is assumed that  $g>0\,$ if $\ds\,\sum^n_{m=1}|y_m|>0$.
Below we will describe some possible ways of choosing the function $g=g(x,y_1,\dots,y_n)$.

\subsection{Special cases of non-local transformations}\label{ss:9.2}

Let us consider some possible ways of choosing the function~$g$ in the system~\eqref{eq:08b}.

\begin{description}

\item{$1^\circ$.}
We can take
\begin{equation}
g=\Bl[c_0+\sum^n_{m=1}c_m|f_m(x,y_1,\dots,y_n)|^s\Br]^{1/s},\ \ c_0>0,\ \ c_m>0,\ \ s>0.
\label{eq:08e}
\end{equation}

\end{description}

\medskip

In particular, if we set $c_0=c_m=s=1$ ($m=1,\dots,n$)  in~\eqref{eq:08e}, then
the system~\eqref{eq:08b} takes the form
\begin{align}
\frac{dx}{d\xi}=\frac{1}{\ds 1{+}\sum^n_{m=1}|f_m(x,y_1,\dots,y_n)|},\ \
\frac{dy_m}{d\xi}=\frac{f_m(x,y_1,\dots,y_n)}{\ds 1{+}\sum^n_{m=1}|f_m(x,y_1,\dots,y_n)|},
\label{eq:08f}
\end{align}
where $m=1,\dots, n$.

Unlike the right-hand sides of the original system~\eqref{eq:08},
the right-hand sides of the system~\eqref{eq:08b} with~\eqref{eq:08e}
have no singularities since all the derivatives
are bounded, $|(y_m)'_\xi|\le 1$ ($m=1,\dots, n$); we recall that for blow-up solutions
at least one of the derivatives $(y_m)'_x$  tends to infinity as $x\to x_*$.

The numerical solution of the problem \eqref{eq:08b}--\eqref{eq:08c} with~\eqref{eq:08e}
can be obtained, for example, applying the Runge--Kutta method or other standard numerical methods,
see above.

\begin{description}

\item{$2^\circ$.}
For the system~\eqref{eq:08b} with
\begin{equation}
g=\Bl[1+\sum^n_{m=1}f^2_m(x,y_1,\dots,y_n)\Br]^{1/2},
\label{eq:08d*}
\end{equation}
we get the method of the arc-length transformation~\cite{mor1979}
(the function~\eqref{eq:08d*} is a particular case of~\eqref{eq:08e}
with $c_0=c_m=1$ and $s=2$).
Therefore for blow-up problems, the method based on introducing the non-local variable~\eqref{eq:08d}
is more general than the method of the arc-length transformation.

\item{$3^\circ$.}
In the general case, it is not known in advance whether the solution of the Cauchy problem
\eqref{eq:08}--\eqref{eq:08a} is a solution with usual properties, or is a blow-up solution.
Therefore, in the first stage, the problem \eqref{eq:08}--\eqref{eq:08a} can be solved by any standard fixed-step numerical
method, for example, by the Runge--Kutta method.
If one of the components, for example,  $y_k$, begins to grow very rapidly
(and increases faster than exponential and faster than the other components),
then a hypothesis arises that the corresponding solution is a blow-up solution.
Numerical confirmation of this hypothesis is a rapid growth of the ratio $|f_k/y_k|$
with increasing of the integration region with respect to~$x$.
In this case, it is reasonable to choose the function~$g$ in~\eqref{eq:08b}, for example, as follows:
\begin{align}
g=\frac 1{y_k}f_k(x,y_1,\dots,y_n).
\label{eq:08g}
\end{align}
As a result, the $(k+1)$-th equation of the system~\eqref{eq:08b} is easily integrated and,
taking into account the corresponding initial condition~\eqref{eq:08c}, we arrive at the dependence
\begin{align}
y_k=y_{k0}e^{\xi}.
\label{eq:08h}
\end{align}
Substituting the relations \eqref{eq:08g} and \eqref{eq:08h} into the remaining equations
of the system~\eqref{eq:08b}, we obtain the Cauchy problem
\begin{equation}
\begin{gathered}
\frac{dx}{d\xi}=\frac{y_k}{f_k(x,y_1,\dots,y_n)},\quad \ \frac{dy_m}{d\xi}=\frac{y_kf_m(x,y_1,\dots,y_n)}{f_k(x,y_1,\dots,y_n)},\\
m=1,\dots, n;\quad \ \text{$m\not=k$} \quad \  (\xi>0),
\end{gathered}
\label{eq:08i}
\end{equation}
with the initial conditions \eqref{eq:08c}.

In the right-hand sides of the system~\eqref{eq:08i}, the function $y_k$
should be replaced by the right-hand side of the formula~\eqref{eq:08h}.

The numerical solution of the system~\eqref{eq:08i} with~\eqref{eq:08h} and the initial conditions~\eqref{eq:08c}
can be obtained, for example, by applying the Runge--Kutta method or the other standard numerical methods
with a sufficiently large stepsize in~$\xi$.

\end{description}
\medskip

\begin{example}
We consider the test Cauchy problem for the system of three equations
\begin{equation}
\begin{gathered}
\frac{dy_1}{dx}=-y_1y_2,\quad \ \frac{dy_2}{dx}=y_2^4y_3,\quad \ \frac{dy_3}{dx}=-2y_1;\\
y_1(0)=y_2(0)=y_3(0)=1.
\end{gathered}
\label{x-30}
\end{equation}
The exact solution of this problem has the form
\begin{equation}
y_1=1-x,\quad \ y_2=\frac 1{1-x},\quad \ y_3=(1-x)^2.
\label{x-31}
\end{equation}

\begin{figure}
\centering
{\includegraphics[scale=0.34]{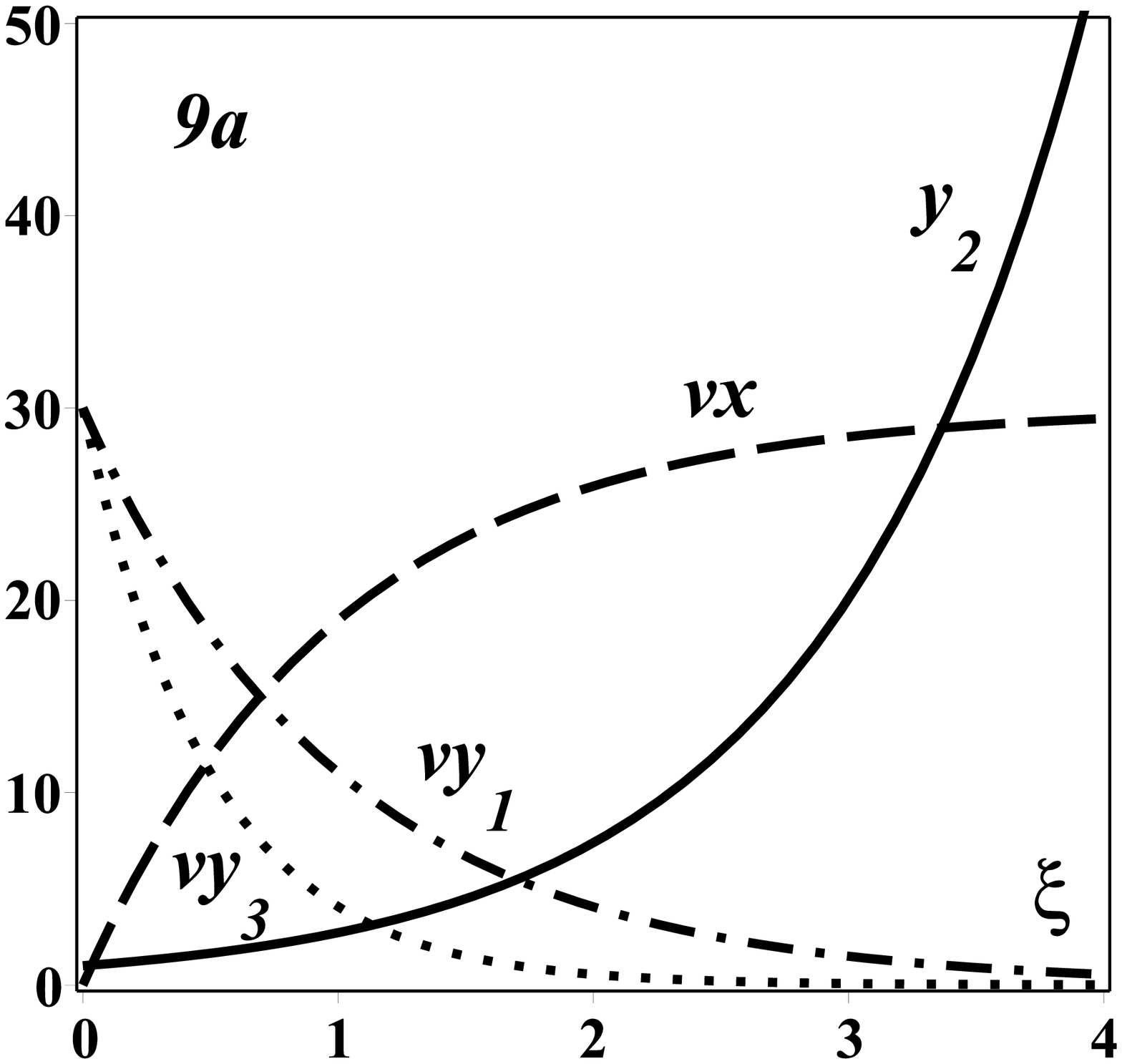}\ \includegraphics[scale=0.34]{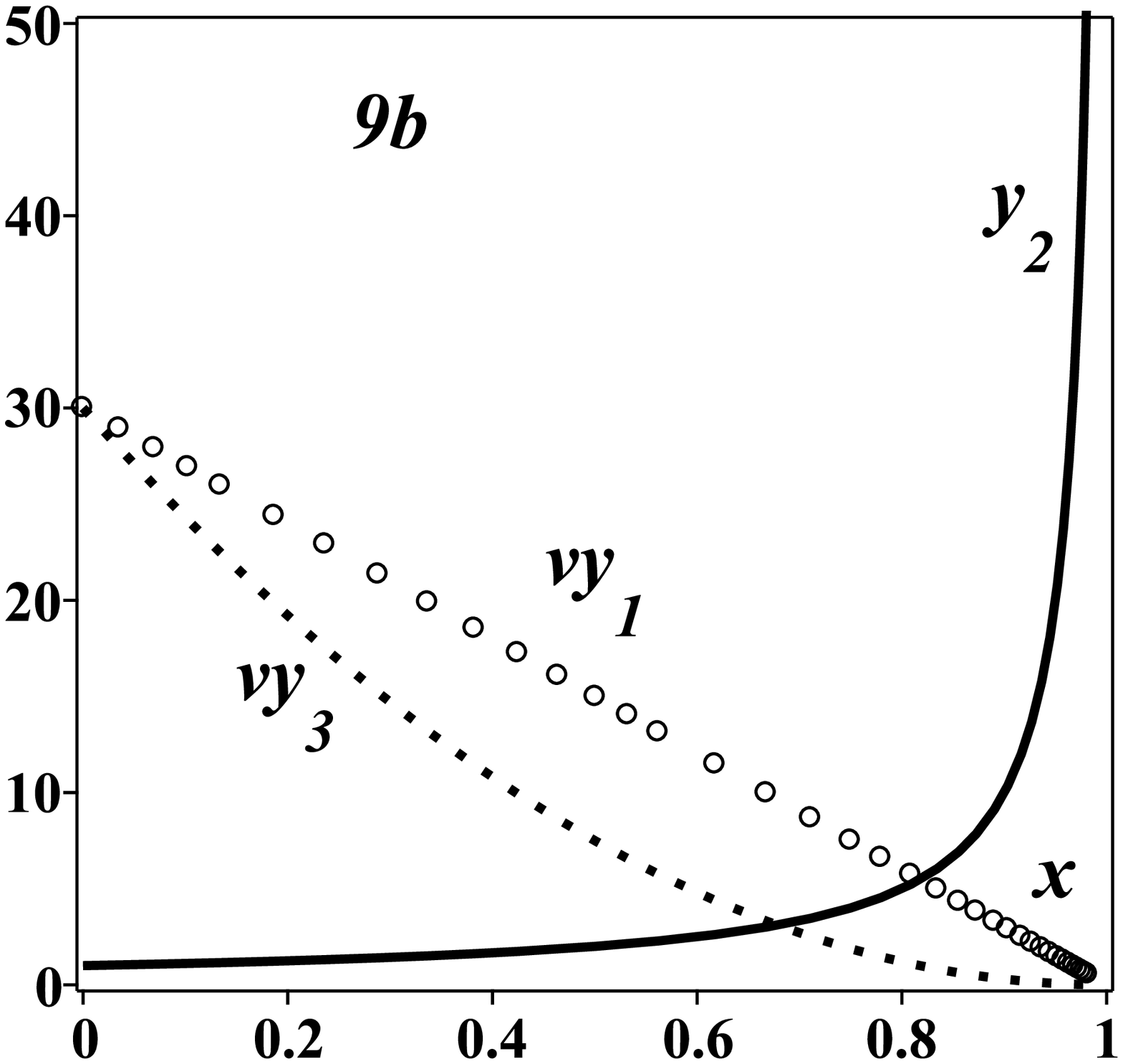}}
\caption{{\it 9a}\dash the dependences $x=x(\xi)$, $y_1=y_1(\xi)$, $y_3=y_3(\xi)$, obtained by numerical solution
of the problem~\eqref{x-32} ($\nu=30$), and $y_2=e^\xi$; {\it 9b}\dash numerical solution of the problem~\eqref{x-32}:
$y_1(x)$ (circles), $y_2(x)$ (solid line), and $y_3(x)$ (points).}
\label{fig:Fig9}
\end{figure}

A trial numerical integration of the problem~\eqref{x-30} by the Runge--Kutta method
shows that the component~$y_2$ grows faster (in magnitude) that the other components.
Using the formulas \eqref{eq:08g} and~\eqref{eq:08h}, we obtain that $g=y_2^3y_3$ and $y_2=e^\xi$.
Substituting these functions into~\eqref{eq:08i}, and taking into account that $f_1=-y_1y_2$ and $f_3=-2y_1$,
we arrive at the equivalent Cauchy problem
\begin{equation}
\begin{gathered}
\frac{dx}{d\xi}=\frac {e^{-3\xi}}{y_3},\quad \ \frac{dy_1}{d\xi}=-\frac {e^{-2\xi}y_1}{y_3},\quad \ \frac{dy_3}{d\xi}=-\frac {2e^{-3\xi}y_1}{y_3};\\
x(0)=0,\quad \ y_1(0)=y_3(0)=1.
\end{gathered}
\label{x-32}
\end{equation}

Unlike the original problem~\eqref{x-30}, the problem~\eqref{x-32} does not have blow-up singularities.
Its exact solution is written in parametric form as follows:
\begin{equation}
x=1-e^{-\xi},\quad y_1=e^{-\xi},\quad \ y_2=e^{\xi},\quad \ y_3=e^{-2\xi}.
\label{x-33}
\end{equation}

The numerical solution of the problem~\eqref{x-32} is shown in Fig.~\ref{fig:Fig9}.
We do not present here the exact dependences~\eqref{x-33}, since they almost coincide
(up to the maximum error~$0.025\%$) with the results of the numerical solution.
\end{example}

\medskip

\begin{remark}
In the methods described in Items~$1^\circ$ and~$2^\circ$, the rate of approximation
of the function $x=x(\xi)$ to the asymptote, that determines the singular point~$x_*$,
will be power-law behavior with respect to~$\xi$, while the method presented in Item~$3^\circ$,
yields the exponential rate of approximation of the singular point.
\end{remark}

\begin{remark}
For systems of equations \eqref{eq:08} of polynomial type, the most growing component $y_k$
can be determined by substituting the approximate functions
$y_1=\alpha_1(x_*-x)^{-\beta_1}$, \dots, $y_n=\alpha_n(x_*-x)^{-\beta_n}$
into the equations.
Then, from the analysis of the obtained algebraic relations, the largest exponent
$\beta_k=\max[\beta_1,\dots,\beta_n]$ is found, where $\beta_k>0$.
The component $y_k$ is used in formula \eqref{eq:08g} for the function~$g$.
\end{remark}

\begin{example}
Consider the problem \eqref{x-30}.
The solution in the neighborhood of the singular point is sought in the form
\begin{equation}
y_1=\alpha_1(x_*-x)^{-\beta_1},\quad \ y_2=\alpha_2(x_*-x)^{-\beta_2},\quad \ y_3=\alpha_3(x_*-x)^{-\beta_3}.
\label{eqeq}
\end{equation}
Substituting the expressions \eqref{eqeq} into \eqref{x-30}, we obtain a simple system of linear algebraic equations
for the exponents $\beta_m$ ($m=1,\,2,\,3$).
The solution of the system is
$$
\beta_1=-1,\quad \ \beta_2=1, \quad \ \beta_3=-2.
$$
The maximum exponent is $\beta_2$. Therefore, the component $y_2$ should be used
for the function~$g$ in formula~\eqref{eq:08g}.
\end{example}

\begin{remark}
If the two components, $y_k$ and $y_j$, simultaneously have a blow-up behavior
(with the same or different rate of approaching to infinity as $x\to x_*$), then we also can choose, for example,
$g=c+|f_k|+|f_j|$ or $g=c_1+\sqrt{c_2+f_k^2+f_j^2}$ in~\eqref{eq:08b}.
Here $c$, $c_1$, and $c_2$ are some non-negative constants.
\end{remark}

\begin{remark}
The technique developed in Section~9 can also be used in
Cauchy problems for partial differential equations (PDEs) with blow-up solutions,
if to apply the methods leading to systems of ODEs (for example, in
projection methods and the method of lines \cite{sch1991,sch2009}).
\end{remark}

\subsection{Method based on differential constraints}\label{ss:9.3}

For the Cauchy problems that are described by the systems of ODEs and have blow-up solutions,
the method of numerical integration, based on introducing a non-local variable,
can be generalized if, instead of the relation~\eqref{eq:08d} to take
the first-order differential constraint
\begin{align}
\xi'_x=g(x,y_1,\dots,y_n,\xi)
\label{eq:xz***}
\end{align}
with the initial condition $\xi(x=x_0)=\xi_0$.

If we set $\xi_0=0$, then the use of the differential constraint~\eqref{eq:xz***}
leads to the problem~\eqref{eq:08b}, where the function $g(x,y_1,\dots,y_n)$
must be replaced by the function $g(x,y_1,\dots,y_n,\xi)$ in the equations.

\medskip

\section{Blow-up problems for higher-order ODEs}\label{s:10}

\subsection{Reduction of higher-order ODEs to a system of first-order ODEs}\label{ss:10.1}

Consider the Cauchy problem for the $n$\,th-order ODE:
\begin{equation}
\begin{gathered}
y^{(n)}_x=f(x,y,y^{\prime}_x,\dots,y^{(n-1)}_x)\quad \ (x>x_0);\\
y(x_0)=y_0,\quad y^{\prime}_x(x_0)=y^{(1)}_0,\quad\ldots,\quad y^{(n-1)}_{x}(x_0)=y^{(n-1)}_0,
\end{gathered}
\label{eq:300}
\end{equation}
where $y^{(k)}_x=d^ky/dx^k$ ($k=3,\dots,n$).

The Cauchy problem for one $n$\,th-order ODE~\eqref{eq:300} is equivalent to
the Cauchy problem for a system of $n$~coupled first-order equations of the special form
\begin{equation}
\begin{gathered}
y^{\prime}_1=y_2,\quad \ y^{\prime}_2=y_3,\quad \ \dots,\quad \ y^{\prime}_{n-1}=y_{n},\quad \ y_{n}'=f(x,y_1,y_2,\ldots,y_{n});\\
y_1(x_0)=y_0,\quad y_2(x_0)=y^{(1)}_0,\quad\ldots,\quad y_n(x_0)=y^{(n-1)}_0,
\end{gathered}
\label{eq:301}
\end{equation}
where the prime denotes the derivative with respect to $x$ and $y_1\equiv y$.

The problem~\eqref{eq:301} is a particular case of the Cauchy problem \eqref{eq:08}--\eqref{eq:08a} and
the general methods described in Sections~9.1 and 9.2 are applicable to it.

\subsection{Blow-up problems for third-order ODEs}\label{ss:10.2}

Let us consider the Cauchy problem for the nonlinear third-order ODE of the general form
\begin{align}
y^{\prime\prime\prime}_{xxx}=f(x,y,y^\prime,y^{\prime\prime}_{xx})\quad (x>0);
\quad y(0)=y_0,\ y^{\prime}_x(0)=y_1,\ y^{\prime\prime}_{xx}(0)=y_2.
\label{ee:001}
\end{align}

The problem for one third-order ODE~\eqref{ee:001} is equivalent to
the following problem for the system of three coupled first-order equations:
\begin{equation}
\begin{gathered}
y^{\prime}_x=t,\quad \ t^\prime_x=w,\quad \ w^\prime_x=f(x,y,t,w)\quad \ \ (x>0);\\
y(0)=y_0,\quad \ t(0)=y_1,\quad \ w(0)=y_2.
\end{gathered}
\label{ee:002}
\end{equation}
The introduction of the non-local variable~\eqref{eq:08d} transforms the system~\eqref{ee:002} to the form
\begin{equation}
\begin{gathered}
x^\prime_\xi=\frac 1{g},\quad \
y^{\prime}_\xi=\frac tg,\quad \ t^\prime_\xi=\frac wg,\quad \ w^\prime_\xi=\frac{f}g\quad \ \ (\xi>0);\\
x(0)=0,\quad \ y(0)=y_0,\quad \ t(0)=y_1,\quad \ w(0)=y_2,
\end{gathered}
\label{ee:003}
\end{equation}
where $f=f(x,y,t,w)$ and $g=g(x,y,t,w)$.

\medskip
Let us consider various possibilities for choosing the function~$g$ in the system~\eqref{ee:003}.

\begin{description}

\item{$1^\circ$.} We can take $g=\bl(c_1+c_2|t|^s+c_3|w|^s+c_4|f|^s\br)^{1/s}$ for $c_m>0$ and $s>0$. The case $c_1=c_2=c_3=c_4=1$ and $s=2$
corresponds to the method of the arc-length transformation~\cite{mor1979}.

\item{$2^\circ$.} We can take $g=t/y$ (or $g=kt/y$, where $k>0$ is a constant).
In this case, the system~\eqref{ee:003} is simplified, since the second equation is directly integrated,
and taking into account the second initial condition, we obtain $y=y_0e^\xi$.

\item{$3^\circ$.} We can take $g=w/t$ (or $g=kw/t$ with $k>0$).
In this case, the system~\eqref{ee:003} is simplified, since the third equation is directly integrated,
and we obtain $t=y_1e^\xi$. Taking into account the relations \eqref{ee:002}, we also have
$$
\xi=\int^x_{x_0}\frac{w}t\,dx=\int^x_{x_0}\frac{y^{\prime\prime}_{xx}}{y^{\prime}_{x}}\,dx=\ln\frac{y^{\prime}_{x}}{y_1}.
$$
Thus, this non-local transformation coincides with the modified differential transformation, which was considered in Section~2.3.

\item{$4^\circ$.} Also, we can take $g=f/w$ (or $g=kf/w$ with $k>0$).
In this case, the system~\eqref{ee:003} is also simplified, since the fourth equation is directly integrated,
and we obtain $w=y_2e^\xi$.

\end{description}

The transformations corresponding to the last three cases, $2^\circ$, $3^\circ$, and $4^\circ$,
will be called the special exp-type transformations, they lead
to the solutions, in which the variable~$x$ tends exponentially rapidly to a blow-up point~$x_*$.

\begin{example}
We consider in more detail the test Cauchy problem of the form
\begin{align}
y^{\prime\prime\prime}_{xxx}=6y^4\quad \ \ (x>0);\quad y(0)=y^{\prime}_x(0)=1,\quad y^{\prime\prime}_{xx}(0)=2.
\label{ee:004}
\end{align}
The exact solution of this problem is determined by the formula~\eqref{eq:02g}.

In Table~3, a comparison of the efficiency of the numerical integration methods,
based on various non-local transformations of the form~\eqref{eq:05} is presented by using the example of the test blow-up problem
for the third-order ODE~\eqref{ee:004}.
The comparison is based on the number of grid points needed to make calculations with the same
maximum error (e.g., equal to $0.1$ and~$0.01$).


\begin{table}[!ht]
\begin{center}
\begin{tabular}{|l|l|r|r|r|l|l|r|r|r|}
\hline
\multicolumn{5}{|c|}{\vphantom{} {\small Error$_{{\rm max}}, \% =0.1$} }     \\ \hline \hline
{\small Transformation}     &  \quad {\small Function $g$}  &  {\small Max. interval}  & {\small Stepsize}      & {\small Grid points} \\
{}  &    &  {\small $\xi_{{\rm max}}\qquad$}  & {\small $h\quad$}  & {\small number}  $N$
\\ \hline
{\small Arc-length, Item~$1^\circ$}        &  {\small $g{=}\sqrt{1{+}t^2{+}w^2{+}f^2}$} & 249600.000  & 0.7800 & 320000 \\ \hline
{\small Nonlocal, Item~$1^\circ$}  &  {\small $g{=}1{+}|t|{+}|w|+|f|$}                  & 252000.000  & 1.4000 & 180000 \\ \hline
{\small Hodograph}  &  {\small $g{=}t$}                                                 &     49.010  & 0.1690 & 290  \\ \hline
{\small Special exp-type, Item~$4^\circ$} &  $g{=}f/w$                                  &     11.741  & 0.2498 & 47  \\ \hline
{\small Special exp-type, Item~$2^\circ$} &  $g{=}t/y$                                  &      3.912  & 0.0978 & 40   \\ \hline
{\small Special exp-type, Item~$3^\circ$} &  $g{=}w/t$                                  &      7.828  & 0.2060 & 38   \\ \hline

\hline
\multicolumn{5}{|c|}{\vphantom{} {\small Error$_{{\rm max}}, \% =0.01$} }     \\ \hline \hline
{\small Transformation}     &  \quad {\small Function $g$}  &  {\small Max. interval}  & {\small Stepsize}      & {\small Grid points} \\
{}  &    &  {\small $\xi_{{\rm max}}\qquad$}  & {\small $h\quad$}  & {\small number}  $N$
\\ \hline
{\small Arc-length, Item~$1^\circ$}        &  {\small $g{=}\sqrt{1{+}t^2{+}w^2{+}f^2}$} & 252000.000  &  0.4500 & 560000  \\ \hline
{\small Nonlocal, Item~$1^\circ$}  &  {\small $g{=}1{+}|t|{+}|w|+|f|$}                  & 253580.000  &  0.8180 & 310000  \\ \hline
{\small Hodograph}  &  {\small $g{=}t$}                                                 &     49.020  &  0.0950 & 516  \\ \hline
{\small Special exp-type, Item~$4^\circ$} &  $g{=}f/w$                                  &     11.738  &  0.1381 & 85  \\ \hline
{\small Special exp-type, Item~$2^\circ$} &  $g{=}t/y$                                  &      3.920  &  0.0560 & 70   \\ \hline
{\small Special exp-type, Item~$3^\circ$} &  $g{=}w/t$                                  &      7.827  &  0.1223 & 64   \\ \hline
\end{tabular}
\caption{Various types of analytical transformations applied for numerical
integration of the problem~\eqref{ee:004} with a given accuracy (percent errors are $0.1$ and $0.01$
for $\Lambda_\text{m}\le 50$) and their basic parameters  (maximum interval,
stepsize, grid points number).}
\label{tab:Transfs-ODE3}
\end{center}
\end{table}






It can be seen that the arc-length transformation is the least effective,
since the use of this transformation is associated with a large number of grid points.
In particular, when using the last three transformations, you need 6580--8750 times less
of a number of grid points.
The hodograph transformation has an intermediate (moderate) efficiency.
The use of the last three special exp-type transformations with $g=f/w$, $g=t/y$, and $g=w/t$
gives rather good results. Note that an even smaller number of grid points can be
obtained by using suitable differential constraints.
\end{example}

\section{Elementary approaches allowing one to find the form of new variables}\label{s:11}

We now describe an elementary approach, based on simple semi-geometric considerations,
which allows us to find the form of new variables that transform the original blow-up problem
to a problem, more suitable for numerical integration, that does not have blow-up singularities.

\subsection{Combination of point transformation and hodograph transformation}\label{ss:11.1}

Let us consider the approximate relation~\eqref{eq:01} as an equation connecting
the variables $x$ and $y$. Solving it with respect to~$x$ (for concreteness, we assume that $A>0$), we obtain
\begin{equation}
x=x_*-B_1y^{-1/\beta},\quad \ \ B_1=A^{1/\beta}.
\label{202}
\end{equation}
It is seen that $x$ tends to the blow-up point~$x_*$ slowly enough as $y\to\infty$
(by the power law $\sim y^{-1/\beta}$).
If we make the substitution $y=e^{\xi}$, then the rate of approximation
to the desired asymptotic value~$x_*$ will become exponential with respect to the new variable~$\xi$
(i.e., will increase significantly).
It is convenient to represent the described procedure in the form of a transformation
\begin{equation}
\xi=\ln y,\quad \ z=x,
\label{203}
\end{equation}
where $z=z(\xi)$ is the new unknown function.
As a result, we arrive at the dependence $z=x_*-B_1e^{-\xi/\beta}$,
which can be written in the parametric form
\begin{equation}
x=x_*-B_1e^{-\xi/\beta},\quad \ y=e^{\xi}.
\label{204}
\end{equation}
The transformation \eqref{203} is a combination of two simple point transformations:
1) the non-linear transformation $\bar x=x$, $\bar y=\ln y$ and 2) the hodograph transformation $\xi=\bar y$, $z=\bar x$.
The transformation~\eqref{203} is equivalent to the transformation~\eqref{*}
if $g=f/y$ (see Section~3.1, Item~$5^\circ$) and to the transformation~\eqref{eq:05} if
$g=t/y$ (see Section~7.1, Item~$5^\circ$).

\subsection{Combination of transformation, based on a differential variable, and point transformation}\label{ss:11.2}

Differentiating the asymptotics~\eqref{eq:01}, we have the following relations:
\begin{equation}
y^{\prime}_x=A\beta(x-x_*)^{-\beta-1},\quad \ y=A\BL(\frac{y^{\prime}_x}{A\beta}\BR)^{\ts\!\frac \beta{\beta+1}}
\label{205}
\end{equation}
(the second relation is obtained from the first one after elimination of~$x$ by means of~\eqref{eq:01}).

Excluding~$y$ from~\eqref{202} with the help of the second relation~\eqref{205}, we obtain
\begin{equation}
x=x_*-B_2(y^{\prime}_x)^{\ts-\frac 1{\beta+1}},\quad \ B_2=(A\beta)^{\ts\frac 1{\beta+1}}.
\label{206}
\end{equation}
It is seen that $x$ tends to the blow-up point~$x_*$ slowly enough as $y^{\prime}_x\to\infty$
(in accordance with the power law $\sim (y^{\prime}_x)^{-1/(\beta+1)}$).

If we make the substitution $y^{\prime}_x=e^{\xi}$, then the rate of approximation
to the desired asymptotic value~$x_*$ will become exponential with respect to the new variable~$\xi$
(i.e., will increase significantly).
The described procedure can be represented as a transformation
\begin{equation}
\xi=\ln y^{\prime}_x,\quad \ x=x(\xi),\quad \ y=y(\xi),
\label{207}
\end{equation}
which is based on a combination of the differential transformation $t=y^{\prime}_x$
(see Sections 2.1 and~6.1) and the point transformation $\xi=\ln t$.

The transformation~\eqref{207} determines the asymptotics of the solution~\eqref{eq:01}
in a neighborhood of the blow-up singularity in the parametric form
\begin{equation}
x=x_*-B_2e^{-\frac 1{\beta+1}\xi},\quad \ y=AB_2^{-\beta}e^{\frac \beta{\beta+1}\xi}.
\label{208}
\end{equation}

In Section 2.3, apart from other considerations, it was described how one can obtain
a transformation of the type~\eqref{207}.

\subsection{Relation allowing one to control the calculation process}\label{ss:11.3}

We now derive a useful formula that makes it possible to control the calculation process.

Taking into account the relations \eqref{eq:01} and~\eqref{205},
we differentiate the relation $y/y^{\prime}_x$. After elementary transformations, we obtain
\begin{equation}
\frac 1{\beta}=\frac{yy^{\prime\prime}_{xx}}{(y^\prime_x)^2}-1
=\frac{y}{y^\prime_\xi}\BL(\frac{y^{\prime\prime}_{\xi\xi}}{y^\prime_\xi}-\frac{x^{\prime\prime}_{\xi\xi}}{x^\prime_\xi}\BR)-1,
\label{201}
\end{equation}
where $x=x(\xi)$, $y=y(\xi)$ is the representation of the solution in the parametric form.

For blow-up problems with a power singularity, the constant~$\beta$ must be greater than zero.
Therefore, for numerical representation of solutions in the parametric form, for large values of~$\xi$
the right-hand side of~\eqref{201} must tend to a positive constant (asymptote), which allows us to control the calculation process.

\section{Brief conclusions}\label{s:12}

Three new methods of numerical integration of Cauchy problems
for nonlinear ODEs of the first and second-order, which have a blow-up solution,
are described. These methods are based on differential and non-local
transformations, and also on differential constraints, that lead to
the equivalent problems for systems of equations, whose solutions
are represented in parametric form and have no singularities.

It is shown that:

\begin{description}

\item{(i)} the method based on a non-local transformation of the general form
includes themselves, as particular cases, the method
based on the hodograph transformation, the method of the arc-length transformation, and
the methods based on the differential and modified differential transformations;

\item{(ii)} the methods based on the exp-type and modified differential transformations are
much more efficient than the method based on the hodograph transformation, the method of the
arc-length transformation, and the method based on the differential transformation;

\item{(iii)} the method based on the differential constraints is the most general of the proposed methods.

\end{description}

In the Cauchy problems described by the first-order equations,
two-sided theoretical estimates are established for the critical value of the independent
variable~$x=x_*$, when an unlimited growth of the solution occurs as approaching it.

It is shown that the method based on a non-local transformation of the general form
as well as the method based on the differential constraints admit generalizations
to the $n$\,th-order ordinary differential equations and systems of coupled differential
equations.

\textbf{Note.}  In the near future, the authors are going to submit 
an article on the topic 
``Numerical integration of non-monotonic blow-up problems based on non-local transformations''.


\begin{thebibliography}{000}
\providecommand{\natexlab}[1]{#1}
\providecommand{\url}[1]{\texttt{#1}}
\expandafter\ifx\csname urlstyle\endcsname\relax
  \providecommand{\doi}[1]{doi: #1}\else
  \providecommand{\doi}{doi: \begingroup \urlstyle{rm}\Url}\fi



\bibitem{stuart1990}
M. Stuart, M. S. Floater.
On the computation of blow-up,
\textit{European J. Applied Math.}, 1990, vol. 1, No. 1, pp. 47--71.

\bibitem{als2005}
E. A. Alshina, N. N. Kalitkin, P.V. Koryakin.
Diagnostics of singularities of exact solutions in computations with error control (in Russian),
\textit{Zh. Vychisl. Mat. Mat. Fiz.}, 2005, vol.~45, No.~10, pp.~1837--1847;
http://eqworld.ipmnet.ru/ru/solutions/interesting/alshina2005.pdf

\bibitem{aco2002}
G. Acosta,  G. Dur\'an, J. D. Rossi.
An adaptive time step procedure for a parabolic problem with blow-up.
\textit{Computing}, 2002, vol. 68, pp. 343--373.

\bibitem{mor1979}
S. Moriguti, C. Okuno, R. Suekane, M. Iri, K. Takeuchi.
\textit{Ikiteiru Suugaku -- Suuri Kougaku no Hatten} (in Japanese).
Baifukan, Tokyo, 1979.

\bibitem{hir2006}
C. Hirota, K. Ozawa.
Numerical method of estimating the blow-up time and rate of the solution of ordinary differential
equations -- An application to the blow-up problems of partial differential equations,
\textit{J. Comput. \& Applied Math.}, 2006, vol. 193, No. 2, pp. 614--637.


\bibitem{meyer98}
R. Meyer-Spasche, Difference schemes of optimum degree of implicitness for a
family of simple ODEs with blow-up solutions, \textit{J. Comput. and Appl.
Mathematics}, 1998, Vol.~97, pp.~137--152.

\bibitem{gor1999}
A. Goriely and C. Hyde, Finite time blow-up in dynamical systems, \textit{Phys. Letters A}, 1999, Vol. 250, pp. 4--6.

\bibitem{eli2006}
U. Elias and H. Gingold, Critical points at infinity and blow-up of solutions of autonomous polynomial differential systems via compactification,
\textit{J. Math. Anal. Appl.,} 2006, Vol. 318, pp. 305--322.

\bibitem{bar2008}
J. Baris, P. Baris, and B. Ruchlewicz, Blow-up solutions of quadratic differential systems, \textit{J. Math. Sciences}, 2008, Vol. 149, No. 4, pp~1369--1375.

\bibitem{nas2009}
N. R. Nassif, N. Makhoul-Karamb, and Y. Soukiassian, Computation of blowing-up solutions for second-order differential equations using re-scaling techniques,
\textit{J. Comput. Applied Math.,} 2009, Vol.~227, 185--195.

\bibitem{dlam2012}
P. G. Dlamini, M. Khumalo.
On the computation of blow-up solutions for semilinear ODEs and parabolic PDEs,
\textit{Math. Problems in Eng.}, 2012, vol. 2012, Article ID 162034, 15~p.

\bibitem{zhou2016}
Y.C. Zhou, Z. W. Yang, H. Y. Zhang, and Y. Wang,
Theoretical analysis for blow-up behaviors of differential equations with piecewise constant arguments, \textit{Appl. Math. Comput.,} 2016,
Vol. 274, pp.~353--361.

\bibitem{tak2017}
A. Takayasu, K. Matsue, T. Sasaki, K. Tanaka,
M. Mizuguchi, S. Oishi. Numerical validation of blow-up solutions of ordinary
differential equations. \textit{J. Comput. Applied Mathematics}, 2017, vol. 314, pp. 10--29.

\bibitem{but}
J. C. Butcher. \textit{The Numerical Analysis of Ordinary Differential Equations:
Runge--Kutta and General Linear Methods}. Wiley-Interscience, New York, 1987.

\bibitem{fox87}
L. Fox and D. F. Mayers, \textit{Numerical Solution of Ordinary Differential Equations
for Scientists and Engineers}, Chapman \& Hall, 1987.

\bibitem{dormand96}
J. D. Lambert, \textit{Numerical Methods for Ordinary Differential Systems\/}, Wiley, New York, 1991.

\bibitem{sch}
W. E. Schiesser. \textit{Computational Mathematics in Engineering and Applied
Science: ODEs, DAEs, and PDEs}.
CRC Press, Boca Raton, 1994.

\bibitem{shampine94}
L. F. Shampine, \textit{Numerical Solution of Ordinary Differential Equations}, Chapman \& Hall/CRC Press, Boca Raton, 1994.

\bibitem{asch}
U. M. Ascher, L. R. Petzold.
\textit{Computer Methods for Ordinary Differential Equations and Differential-Algebraic Equations}.
SIAM, Philadelphia, 1998.

\bibitem{korn}
G. A. Korn,  T. M. Korn.
\textit{Mathematical Handbook for Scientists and Engineers, 2nd Edition}.
New York: Dover Publ., 2000.

\bibitem{shing2009}
I. K. Shingareva, C. Liz\'arraga-Celaya.
\textit{Maple and Mathematica. A Problem Solving Approach for Mathematics, 2nd Edition.}
Springer, Wien -- New York, 2009.

\bibitem{grif2010}
\textit{D. Griffiths, D. J. Higham}. Numerical Methods for Ordinary Differential Equations.
Springer, Wien -- New York, 2010.

\bibitem{zai1993}
V. F. Zaitsev, A. D. Polyanin.
\textit{Handbook of Nonlinear Differential Equations: Applications in Mechanics, Exact Solutions} (in~Russian), Nauka, Moscow, 1993.

\bibitem{zai1994}
V. F. Zaitsev, A. D. Polyanin.
\textit{Discrete-Group Methods for Integrating Equations of Nonlinear Mechanics}.
CRC Press -- Begell House, Boca Raton, 1994.

\bibitem{mur2010}
C. Muriel and J. L. Romero. Nonlocal transformations and linearization of second-order ordinary differential equations,
\textit{J. Physics A: Math. Theor.,} Vol. 43, No. 43, 434025.

\bibitem{kudr2016a}
N. A. Kudryashov and D. I. Sinelshchikov. On the integrability conditions for a
family of Li\'enard-type equations, \textit{Regular and Chaotic Dynamics}, 2016, Vol. 21, No. 5, pp. 548--555.

\bibitem{kudr2016b}
N. A. Kudryashov and D. I. Sinelshchikov. On the criteria for integrability of the Li\'enard equation,
\textit{Applied Mathematics Letters}, 2016, Vol. 57, pp. 114--120.

\bibitem{pol2012}
A.D. Polyanin, V.F. Zaitsev,
Handbook of Nonlinear Partial Differential
Equations, 2nd Edition, Chapman \& Hall/CRC Press, Boca Raton--London, 2012.

\bibitem{mikh1965}
S. G. Mikhlin, Kh. L. Smolitsky. \textit{Approximate methods for solving differential and integral equations}.
(in~Russian), Nauka, Moscow, 1965.

\bibitem{kam}
E. Kamke.
\textit{Differentialgleichungen:  L\"osungsmethoden und  L\"osungen, I, Gew\"ohnliche Differentialgleichungen.}
B.G. Teubner, Leipzig, 1977.

\bibitem{pol2003}
A. D. Polyanin, V. F. Zaitsev.
\textit{Handbook of Exact Solutions for Ordinary Differential Equations, 2nd Edition.}
Chapman \& Hall/CRC Press, Boca Raton -- London, 2003.

\bibitem{mur}
G. M. Murphy.
\textit{Ordinary Differential Equations and Their Solutions.}
D. Van Nostrand, New York, 1960.

\bibitem{kudr}
N. A. Kudryashov.
\textit{Analytical Theory of Nonlinear Differential Equations} (in Russian).
Institute of Computer Science, Moscow -- Izhevsk, 2004.

\bibitem{pol2016}
A. D. Polyanin and A. I. Zhurov.
Parametrically defined nonlinear differential equations and their solutions: Applications in fluid dynamics.
\textit{Applied Mathematics Letters}, 2016, vol.~55, pp.~72--80.

\bibitem{pol2017}
A. D. Polyanin and A. I. Zhurov.
Parametrically defined nonlinear differential equations, differential–algebraic equations, and implicit ODEs:
Transformations, general solutions, and integration methods,
\textit{Applied Mathematics Letters}, 2017, vol. 64, pp. 59--66.

\bibitem{pol_shi2017}
A. D.~Polyanin and I. K. Shingareva.
The use of differential and non-local transformations for numerical integration of non-linear blow-up problems,
\textit{Int. J. Non-Linear Mechanics}, 2017, Vol.~95, pp.~178--184.

\bibitem{sch1991}
W. E. Schiesser, The Numerical Method of Lines, Academic Press, 1991
[see, also, S. Hamdi, W. E. Schiesser, and G. W. Griffiths, Method of lines, Scholarpedia, 2007,  2(7):2859; doi:10.4249/scholarpedia.2859].

\bibitem{sch2009}
W. E. Schiesser and G. W. Griffiths, A Compendium of Partial Differential Equation Models: Method of Lines Analysis with Matlab, Cambridge University Press, 2009.



\end{thebibliography}
\end{document}